\documentclass{ws-jktr2}
\usepackage{graphicx}

\newcommand{\bysame}{\mbox{\rule{3em}{.4pt}}\,}

\usepackage{amssymb}
\usepackage{amsmath}
\usepackage{latexsym}

\def\draftnote{ }
\renewcommand\today{}

\newcommand\ZZ{{\Bbb Z}}

\newcommand\CC{{\Bbb C}}
\newcommand\QQ{{\Bbb Q}}
\newcommand\qed{\hfill $\Box$ \hfill \\} 
\renewcommand\qed{\hfill $\Box$ \hfill \\}

\newcommand\xmtx{
\Bigl[\begin{array}{cc}
1&1\\
0&1
\end{array}\Bigr]}
\newcommand\mtx[4]{
\Bigl[\begin{array}{cc}
#1&#2\\
#3&#4
\end{array}\Bigr]}

\newtheorem{thm}{Theorem}[section]

\newtheorem{lemm}[thm]{Lemma}

\newtheorem{prop}[thm]{Proposition}

\newtheorem{yosou}[thm]{Conjecture}
 \newtheorem{dfn}[thm]{Definition}
\newtheorem{rem}[thm]{Remark}
 \newtheorem{ex}[thm]{Example}

\title{
Evaluations of the twisted Alexander polynomials
of $2$-bridge knots at $\pm 1$}

\author{Mikami Hirasawa}
\address{Department of Mathematics,
Nagoya Institute of Technology,\\
Nagoya Aichi 466-8555 Japan\\
{\it E-mail: hirasawa.mikami@nitech.ac.jp}
}

\author{Kunio Murasugi}
\address{Department of Mathematics,
University of Toronto,\\
Toronto, ON M5S2E4 Canada\\
{\it E-mail: murasugi@math.toronto.edu}
}

\begin{document}

\maketitle
 \pagestyle{myheadings}
 \markboth{
Twisted Alexander polynomials at $\pm 1$
}
 {M. Hirasawa \& K. Murasugi}

\begin{abstract}
Let $H(p)$ be the set of $2$-bridge knots 
$K(r), 0<r<1$,
such that the group $G(K(r))$ of $K(r)$
is mapped onto a non-trivial
free product, $\ZZ/2 * \ZZ/p$, $p$ being odd.
Then there is an algebraic integer $s_0$ 
such that for any 
$K(r)$ in $H(p)$, $G(K(r))$
has a parabolic representation $\rho$ into
$SL(2, \ZZ[s_0]) \subset SL(2,\CC)$.
Let $\widetilde{\Delta}_{\rho, K(r)}(t)$ 
be the twisted Alexander
polynomial associated to $\rho$.
Then we prove that for any $K(r)$ in $H(p)$,
$\widetilde{\Delta}_{\rho,K(r)}(1)=-2s_0^{-1}$ and 
$\widetilde{\Delta}_{\rho, K(r)}(-1)=-2s_0^{-1}\mu^2$,
where $s_0^{-1}, \mu \in \ZZ[s_0]$.
The number $\mu$ can be recursively evaluated.

\medskip

Keywords:
 Alexander polynomial, 
 $2$-bridge knot,
 knot group, 
 parabolic representaion,
 twisted Alexander polynomial, 
 continued fraction.
\end{abstract}

\section{
Introduction and statement of the main theorem}

The twisted Alexander polynomial of a knot $K$
is a significant 
generalization of the classical
Alexander polynomial of 
$K$ \cite{Lin} and so far, many attempts 
have been made to prove that both polynomials
share certain important properties
\cite{KL}, \cite{KL2}, \cite{Cha}, \cite{GKM}, 
\cite{HLN}.  
However, such a generalization is by no means straightforward.  
In fact, there are only few studies on 
the corresponding question to 
one of the fundamental properties of 
the Alexander polynomial : $\Delta_K(1) = 1$ \cite{Liv2}.
In this paper, we give some information on the twisted 
Alexander polynomials of $2$-bridge knots evaluated at 
$t=1$ and $-1$.  
To be more precise, given an odd integer $p$, let $K(r)\ r\in \QQ, 0 < r < 1$, be a $2$-bridge knot 
such that $G(K(r))$, the group of $K(r)$, 
is mapped onto a non-trivial free product, 
$\ZZ/2 * \ZZ/p$ 
and $H(p)$ the set of all $2$-bridge knots with this property.
Then there is an algebraic integer $s_0$ 
such that the group of each knot $K(r)$ in $H(p)$ has a parabolic representation $\rho$ in $SL(2,\ZZ[s_0]) \subset
SL(2,\CC)$ defined by 
$\rho :  x \mapsto  
\Bigl[\begin{array}{cc}
1&1\\
0&1
\end{array}\Bigr]$    and 
$y \mapsto
\Bigl[\begin{array}{cc}
1&0\\
s_0&1
\end{array}\Bigr]$
where $x$ and $y$ are meridian generators of $G(K(r))$.  
Let $\widetilde{\Delta}_{\rho,K(r)}(t)$ be the twisted Alexander polynomial of  $K(r)$ associated to $\rho$. Then we prove: \\

\noindent
{\bf Theorem A.}
{\it  
 For any knot $K(r)$ in $H(p)$, we have:\\
    (1)   $\widetilde{\Delta}_{\rho,K(r)}(1)= -2 s_0^{-1}$, and\\
    (2)   $\widetilde{\Delta}_{\rho,K(r)}(-1)= -2 s_0^{-1}\mu^2$,\\
where both $s_0^{-1}$ and $\mu$ are elements of $\ZZ[s_0]$.
}\\

In particular, $K(1/p)$ belongs to $H(p)$ and since for 
any knot $K(r)$ in $H(p)$, 
there is an epimorphism from $G(K(r))$ in $H(p)$ to $G(K(\frac{1}{p}))$,
it follows that 
$\widetilde{\Delta}_{\rho,K(1/p)}(t)$ divides 
$\widetilde{\Delta}_{\rho,K(r)}(t)$ 
(\cite{KSW} or see Proposition \ref{prop:3.2}(2)),
and the quotient   
$\lambda_{\rho,K(r)}(t) = \widetilde{\Delta}_{\rho,K(r)}(t) 
/ \widetilde{\Delta}_{\rho,K(1/p)}(t)$ 
is a symmetric polynomial over $\ZZ[s_0]$ 
(Proposition \ref{prop:3.2} (3)). 
Then Theorem A, Proposition \ref{prop:2.1} 
and (4.3)(2) 
imply that
$\lambda_{\rho,K(r)}(1) = 1$ and  
$\lambda_{\rho,K(r)}(-1) = \mu^2$, for $\mu \in \ZZ[s_0]$.
If $p = 3$, then $s_0= -1$, and hence   
$\lambda_{\rho,K(r)}(t)$ is the Alexander polynomial 
$\Delta_K(t)$ of some knot K.
However, the second condition of Theorem A
gives a strong restriction 
for $\Delta_K(t)$.  
Therefore, for example, the quadratic Alexander polynomial cannot 
be realized as the polynomial 
$\lambda_{\rho,K(r)}(t)$ for any knot $K(r)$, since
the degree of $\lambda_{\rho, K(r)}(t)$
must be a multiple of $4$ (Proposition \ref{prop:3.3}). 
On the other hand, for some particular $r$,   
$\lambda_{\rho,K(r)}(t)$ can be realized as the Alexander polynomial.  
In fact, we can prove:\\

\noindent
{\bf Proposition 3.5.}
{\it
 For any odd integers $p$ and $q$,
$\lambda_{\rho,K(1/pq)}(t) = \Delta_{K(1/q)}(t^{2p})$. 
}\\

The number  $\mu  \in \ZZ[s_p]$ is a knot invariant, 
and $\mu$ can easily be  evaluated by using a recursion formula.  
(See Proposition \ref{prop:9.1}.)

After the first draft of the present paper 
was completed, 
we learned that D. Silver and 
S. Williams have been 
studying a similar problem with a different motivation
and they propose a quite interesting conjecture 
that is closely related to 
Theorem A.  
As an application of Theorem A, 
we prove their conjecture partially
for 2-bridge knots in $H(p)$ in Section 10.  

This paper is organized as follows.
In Section 2, 
first we give a quick review of the definition of 
the twisted Alexander polynomial 
and state their basic properties.  
Then we define a parabolic representation of a 
$2$-bridge knot $K(r)$ 
and for a few values of $r$, 
we calculate the twisted Alexander 
polynomial of $K(r)$ 
associated to this representation.   
In Section 3, we introduce a polynomial 
$\lambda_{\rho,K(r)}(t)$  
for $K(r)$ when $G(K(r))$ is mapped onto the free product 
$\ZZ/2 * \ZZ/p, p$ being odd, and determine  
$\lambda_{\rho,K(r)}(t)$ for some 
values $r$.  
In Sections 4, we introduce a
$\ZZ[s_0]$-algebra $\widetilde{A}(s_0)$ 
that is our fundamental tool to prove Theorem A, 
and verify
two technical lemmas about $\widetilde{A}(s_0)$.  
In Section 5, as the first step toward the proof of 
Theorem A, 
we show that Theorem A is reduced to two
formulas in the algebra $\widetilde{A}(s_0)$.  
The purpose of the next section, Section 6, is to show 
that we only need to prove Theorem A 
for much restricted rationals $r$. 
(See Propositions \ref{prop:6.3} and \ref{prop:8.1}.)
In Section 7, 
we prove the first part of 
Theorem A, and the second part of Theorem A is proved 
in Section 8. 
In Section 9, we provide an algorithm to evaluate 
the number $\mu$ appeared in Theorem A.  
In the last section, 
Section 10, 
we state Silver-Williams Conjecture and prove 
their conjecture for torus knots $K(1/p), p$ odd, 
and for 2-bridge knots in $H(p)$. 
In Appendix, 
we give an outline of the proofs of 
Proposition  \ref{prop:2.1} 
and  (10.4)(2),
and also give
a proof of Proposition \ref{prop:3.4}.  

\section{Definition and Examples}

In this section, first we quickly review the definition of 
the twisted Alexander polynomials and their properties 
that we will use throughout this paper.  
 For the details, we refer to \cite{W}. 
Later in this section, we define a parabolic representation of 
the group of a $2$-bridge knot $K(r)$.
(See \cite{R2}.)

Let  $\rho : G = G(K) \rightarrow GL(n, \CC)$ 
be a linear representation of the group of  a knot $K$.
Let $G =\langle
x_1, x_2, \cdots, x_m | r_1, r_2, \cdots, r_{m-1}\rangle$ 
be a Wirtinger presentation of $G(K)$.
Denote by $M_{p,q}(R)$ the ring of $p\times q$  
matrices over a ring $R$.
Let $A = 
\left[ \frac{\partial r_i}{\partial x_j}\right]
\in M_{m-1,m} (\ZZ[x_1^{\pm 1}, \cdots , x_m^{\pm 1}])$ 
be the Alexander matrix, where  
$\frac{\partial}{\partial x_j}$ denotes Fox free derivatives and 
$\ZZ[x_1^{\pm 1}, \cdots , x_m^{\pm 1}]$ 
is a non-commutative ring of Laurent polynomials.                                
The
square matrix
$\widehat{A} = 
\left[\frac{\partial r_i}{\partial x_j}\right]_{1\leq i, j \leq {m-1}}$
is obtained by deleting  the last column of $A$.
We define a homomorphism
$\Phi$ from the group ring $\ZZ{G}$ into $M_{n,n}(\CC[t^{\pm 1}])$
by $\Phi(x_i)=\rho(x_i)t$.
Then           
$( \frac{\partial r_i}{\partial x_j})^\Phi \in 
M_{n,n}(\CC[t^{\pm 1}])$,
and hence
$(\widehat{A})^\Phi = 
\left[(\frac{\partial r_i}{\partial x_j})^\Phi\right] \in 
M_{(m-1)n,(m-1)n}(\CC[t^{\pm 1}])$. 

\begin{dfn} \cite{W}\label{dfn:2.1}
The twisted Alexander polynomial of $K$ associated 
to  $\rho$ is defined as follows:
$$\widetilde{\Delta}_{\rho, K}(t) =  
\frac{\det \widehat{A}^\Phi}{\det (x_{m}^ \Phi - 1)}
\in \CC[t^{\pm 1}].$$
If $\rho$ is unimodular,
this is an invariant of $K$ up to $\pm t^{nk}$.   
\end{dfn}

We should note that for any linear representation
$\rho$, the ambiguity of this invariant 
is completely eliminated by Kitayama.  
For the precise formulation, see \cite{Ktym}.  

 \begin{rem}\label{rem:2.1}
(1)  If  $\rho : x_i  \mapsto   I \in GL(n, \CC)$ 
is a trivial representation, then             
$\widetilde{\Delta}_{\rho, K}(t) =  
\left[\frac{\Delta_K(t)}{t - 1}\right]^n$,  
where $\Delta_K(t)$ is  
the Alexander polynomial of a knot $K$.
(2)  In general, $\widetilde{\Delta}_{\rho,K}(t)$ 
is a rational function, but it is shown \cite{W} 
that if the commutator subgroup $G'$ contains 
an element $w$ such that $1$ is not an eigenvalue of  
$\rho (w)$, then $\widetilde{\Delta}_{\rho,K}(t)$ is 
a Laurent polynomial over $\CC$, namely,
 $\widetilde{\Delta}_{\rho,K}(t)  
 \in \CC[t^{\pm 1}]$.
(3) For any presentation $\rho: G(K)\rightarrow
SL(n, \CC)$,
$\widetilde{\Delta}_{\rho,K}(t)$ is symmetric 
\cite{Ki}.
\end{rem}

Now we study parabolic representations of the
$2$-bridge knot groups (c.f. \cite{R2}).
Let $r$ be a rational number, $0<r = \frac{\beta}{\alpha} < 1$, 
where both $\alpha$ and $\beta$ are odd and  
$\gcd(\alpha,\beta)=1$, 
and $K(r)$ is the $2$-bridge knot of type $(\alpha,\beta)$.  

Let $F(x,y)$ be the free group freely generated by $x$ and $y$.
For $k = 1,2, \cdots, \alpha - 1$, let  
$\eta_k = [\frac{k\beta}{\alpha}]$, where  $[\cdot  ]$  
denotes Gaussian symbol and let 
$\varepsilon_k  = (-1)^{\eta_k}$.      

Using the word $W$ given by 

\begin{equation}
 W = x^{\epsilon_1}y^{\epsilon_2} x^{\epsilon_3} y^{\epsilon_4} \cdots  x^{\epsilon_{\alpha -2}} y^{\epsilon_{\alpha -1}},
 \end{equation}
 
we obtain a Wirtinger presentation of $G(K(r))$:

\begin{equation}
G(K(r)) = \langle
x,y| WxW^{-1}y^{-1} = 1\rangle.
\end{equation}
 

For each $r, 0 < r < 1$, 
there is a non-commutative representation \\
$\rho:G(K(r)) \rightarrow SL(2,{\Bbb C})$ 
such that  

\begin{equation}
\rho(x)=   
\Bigl[\begin{array}{cc}
1&1\\
0&1
\end{array}\Bigr]\ 
{\rm and}\ 
\rho(y)=  
\Bigl[\begin{array}{cc}
1&0\\
s_r&1
\end{array}\Bigr],
s_r  \neq 0.
\end{equation}


Here a complex number $s_r$ is determined as follows
\cite{R2}.
Let  $G = G(K(r)) =\langle
x ,y | WxW^{-1}y^{-1} =1\rangle$ 
be a Wirtinger presentation of $G$ given by (2.2).\\
Set  $\rho(x)=
\Bigl[\begin{array}{cc}
1&1\\
0&1
\end{array}\Bigr]$ and  
$\rho (y) = 
\Bigl[\begin{array}{cc}
1&0\\
z&1
\end{array}\Bigr]$,                 
where  $z$ is a variable.  \\
 Compute   
 $\rho (W)
  = 
\Bigl[\begin{array}{cc}
a(z)&b(z)\\
c(z)&d(z)
\end{array}\Bigr]$,
where $a,b,c$ and $d$ are polynomials on $z$.
Then equality $Wx=yW$  yields 
$\Bigl[\begin{array}{cc}
a&a+b\\
c&c+d
\end{array}\Bigr]=
\Bigl[\begin{array}{cc}
a&b\\
za+c&zb+d
\end{array}\Bigr]$. \\
The number $s_r$ we sought is a root of $a(z)=0$
\cite[Theorem 2]{R2}.   
For convenience, we call $\rho$ a 
{\it canonical representation} of $G(K(r))$,
and  $a(z)$ the {\it representation polynomial} of $\rho$.
Since $G' \ni xy^{-1}$ and $\rho (xy^{-1})=
\Bigl[\begin{array}{cc}
1-s_r&1\\
-s_r&1
\end{array}\Bigr]$,
$1$ is not an eigenvalue of  
$\rho (x y^{-1})$, and by Remark \ref{rem:2.1} (2) and (3), 
we see that 
$\widetilde{\Delta}_{\rho,K}(t)$ is a symmetric 
Laurent polynomial over $\ZZ[s_r]$.
It is known \cite{R2} that the representation
polynomial
$a(z)$ is a separable polynomial of degree $\frac{\alpha-1}{2}$.

If $r = 1/p, p = 2n + 1$. 
Then W = $(xy)^n$, and it is easy to show that 
the representation polynomial
$a_n(z)$ is a monic polynomial of degree $n$ and further, 
the constant term is also $1$.     
We study $a_n(z)$ in Section 10.

\begin{ex}\label{ex:2.1}
(1)   Let  $r=1/3$.  
Then  $W=xy$ and hence  $s_r = -1$. 
Therefore, 
$\rho :   x \mapsto 
\Bigl[\begin{array}{cc}
1&1\\
0&1
\end{array}\Bigr]$
and  $y  \mapsto 
\Bigl[\begin{array}{cc}
1&0\\
-1&1
\end{array}\Bigr]$                              
gives a parabolic representation 
$\rho : G(K(1/3)) \rightarrow  SL(2,\ZZ)$.  
A simple computation shows that 
$\widetilde{\Delta}_{\rho, K(1/3)}(t) = 1+t^2$.
 
(2)  Let $r=3/5$.  Then  $W=xy^{-1}x^{-1}y$, 
and hence $s_r = -w$, where $w$ is a primitive cubic root of $1$. 
Thus we have a parabolic 
representation $\rho:  G\rightarrow 
SL(2, \ZZ[w])
\subset SL(2,\CC)$ 
and the twisted Alexander polynomial of $K(3/5)$ 
associated to 
$\rho$ is $\widetilde{\Delta}_{\rho, K(3/5)}(t) = 1-4t+t^2$.

(3) Let  $r = 3/7$, then $W=xyx^{-1}y^{-1}xy$  and 
$s_r$ is a root of  $1+2z +z^2 +z^3 =0$.  
The twisted Alexander polynomial associated 
to this representation is
$\widetilde{\Delta}_{\rho,K(3/7)}(t) = 
-(4+s_{r}^2)+4t-(4+s_{r}^2) t^2$.
\end{ex}

\begin{prop}\label{prop:2.1}
The twisted Alexander polynomial of $K(1/p)$, $p=2n+1$,
associated to a canonical representation $\rho$
is given by 
\begin{equation}
\widetilde{\Delta}_{\rho,K(1/p)}(t)
= 
b_1 + b_2 t^2 + b_3 t^4  + \cdots + b_n t^{2n-2} + 
b_{n} t^{2n}+ b_{n-1} t^{2n+2} +   \cdots + b_1t^{4n-2},
\nonumber
\end{equation}
where $b_k$ is the $(1,2)$-entry of  
$\rho (xy)^k$,
and ${\displaystyle 
b_k=\sum_{j=0}^{k-1}
\binom{ktj}{2j+1}s^{j}_{1/p}}$
(see \cite{S}).
\end{prop}

For a proof, see Appendix (I).

\section{
Twisted Alexander polynomials of $2$-bridge knots
}

Suppose that there is an epimorphism from
$G(K(r))$
to non-trivial free product, $\ZZ/2*\ZZ/p$ for some odd $p$.  
Let $H(p)$ be the set of these knots $K(r)$.  
The following proposition is proved in \cite{GonR}.

\begin{prop}\label{prop:3.1}
Let $K(r)$ be an element of $H(p)$. 
We may assume without loss of generality that
$0<r=\dfrac{\beta}{\alpha}<1$, where
$0<\beta<\alpha$, $\alpha \equiv \beta$ 
$\equiv~1\ ({\rm mod\ } 2)$
and $\gcd(\alpha, \beta)$
$=1$.
Then the continued fraction of $r$ is of the form:\\
$ r = [pk_1, 2m_1, pk_2, 2m_2,\dots, 2m_{q}, pk_{q+1}]$,
where $m_i$ and $k_j$ are non-zero integers.
\end{prop}

Here, the continued fraction of $r$ is defined as follows:
%
%


\medskip
$\ \ r=\dfrac{\beta}{\alpha}=
\dfrac{1}{pk_1-\dfrac{1}{2m_1-
\dfrac{1}{pk_2-
\genfrac{.}{.}{0pt}{}{}{
\genfrac{.}{.}{0pt}{}{}{
\ddots 
\genfrac{.}{.}{0pt}{}{}{
\mbox{ $-$}}
}
\genfrac{.}{.}{0pt}{}{}{
\genfrac{.}{.}{0pt}{}{}{
\dfrac{1}{2m_q -
\dfrac{1}{pk_{q+1}}}}
}
}
}}}$
\medskip

A different characterization of continued fractions of 
$r$ for $K(r)$ in $H(p)$ is given in 
Appendix (IV).

According to \cite{ORS}, 
there is an epimorphism $\varphi$ 
from  $G(K(r)), K(r)\in H(p)$,
onto
$G(K(1/p))$ sending meridians of $K(r)$ to those of $K(1/p)$.  
Therefore, the canonical parabolic representation
$\rho: G(K(1/p))
\rightarrow  SL(2,\ZZ[s_{1/p}]) \subset SL(2, \CC)$ 
defined by

\hfill $\rho(x)= 
\xmtx\ 
{\rm and}\ 
\rho(y)=   
\mtx{1}{0}{s_{1/p}}{1}$\hfill(3.1)

\noindent                           
can be extended  to a parabolic representation

\hfill $\rho\varphi:G(K(r)) \rightarrow G(K(1/p)) 
\rightarrow SL(2,\ZZ[s_{1/p}]) \subset SL(2,\CC)$
\hfill(3.2)

\setcounter{equation}{2}

\noindent  
and we can define the twisted Alexander polynomials of $K(r)$ and 
$K(1/p)$ associated to $\rho\varphi$ and $\rho$, respectively. \\

First we prove the following;

\begin{prop}\cite{KSW} \label{prop:3.2}
Let $\rho\varphi : G(K(r)) \rightarrow  
SL(2,\ZZ[s_{1/p}]) \subset SL(2, \CC)$ be the parabolic representation defined 
by (3.2).  
Then,\\
(1) Both $\widetilde{\Delta}_{\rho,K(1/p)}(t)$ and 
$\widetilde{\Delta}_{\rho\varphi,K(r)}(t)$ are polynomials over 
$\ZZ[s_{1/p}]$, and \\
(2)  $\widetilde{\Delta}_{\rho,K(1/p)}(t)$ divides 
$\widetilde{\Delta}_{\rho\varphi,K(r)}(t)$.    \\
Write $\widetilde{\Delta}_{\rho\varphi,K(r)}(t) =
\lambda_{\rho,K(r)}(t)\widetilde{\Delta}_{\rho, K(1/p)}(t)$.  
Then,\\
(3)   $\lambda_{\rho,K(r)}(t)$ is a symmetric polynomial over 
$\ZZ[s_{1/p}]$, and $\lambda_{\rho, K(r)}(t)$ is unique

 up to
$t^{2k}$.
\end{prop}

\noindent
{\it Proof.}
First, (1) follows from 
Remark \ref{rem:2.1}.
To prove (2), consider Wirtinger presentations  
$G(K(1/p)) = \langle x,y| R_0\rangle$ and
$G(K(r)) = \langle x,y | R\rangle$.  
Since an epimorphism $\varphi$ sends $x$ to $x$ and 
$y$ to $y$, it follows that $R  = 1$ in $G(K(1/p))$.  
Therefore, $R$ is written freely as a product of 
conjugates of $R_0$  and

\begin{equation}
R \equiv \prod_{k=1}^m u_j R^{\epsilon_j}_{0} u_j^{-1},
\end{equation}
where $u_j \in F(x,y)$ and $\epsilon_j = \pm 1$, and $A \equiv B$ means that $AB^{-1}$ is equal to the identity of the free group $F(x,y)$.  Therefore,
$\Phi (\frac{\partial R}{\partial x})$ = 
$\sum_{j=1}^m \epsilon_j u_j^\Phi 
(\frac{\partial R_0}{\partial x})^\Phi$,
where  $\Phi :  \ZZ F(x,y)  \rightarrow
M_{2,2}(\ZZ[s_{1/p}][t^{\pm 1}])$.\\
Now  
\begin{eqnarray*}
\widetilde{\Delta}_{\rho,K(r)}(t)
&=&
\det (\frac{\partial R}{\partial x}^\Phi) / \det (y^\Phi - I)
\\
 &=&
\det\Bigl[\sum_{j=1}^m \epsilon_j u_j^\Phi\Bigr]
\det(\frac{\partial R_0}{\partial x})^\Phi/\det (y^\Phi-I)
\\
 &=&
\det\Bigl[\sum_{j=1}^m \epsilon_j u_j^\Phi\Bigr]
\Bigl[\det (\frac{\partial R_0}{\partial x})^\Phi/
\det (y^\Phi-I)\Bigr]
\\
 &=&
\det \Bigl[\sum_{j=1}^m \epsilon_j u_j^\Phi\Bigr] 
\widetilde{\Delta}_{\rho,K(1/p)}(t).
\end{eqnarray*}

%
%
%
%
%
%
%

This proves (2), and further, we see that 

\begin{equation}
\lambda_{\rho,K(r)}(t)=\det \left[\sum_{j=1}^m \epsilon_j u_j^\Phi\right]. 
\end{equation}

By Remark \ref{rem:2.1},
$\lambda_{\rho,K(r)}(t)$  is a 
symmetric polynomial over $\ZZ[s_{1/p}]$.
This proves (3).
\qed

\begin{ex}\label{ex:3.3}
(1) From Proposition \ref{prop:3.1},
 we see that there are epimorphisms from
$G(K(19/45))$ and $G(K(37/213))$ onto $G(K(1/3))$.  
Straightforward calculations show that
$\lambda_{\rho,K(19/45)}(t)= 
25-72 t + 95 t^2 - 72 t^3 + 25 t^4$   and 
$\lambda_{\rho,K(37/213)}(t) = 
4-16 t  + 28 t^2 - 32 t^3  + 28 t^4  -  
16 t^5 + 8 t^6  - 8 t^7  + 4 t^8  -  8 t^{10}
+ 16 t^{11} - 15 t^{12}  + 16 t^{13}  - 8 t^{14}  + 
4 t^{16}   - 8 t^{17} + 8 t^{18}- 16 t^{19}
+28 t^{20}- 32 t^{21}+28 t^{22}- 16 t^{23}+ 4 t^{24}$.
\end{ex}

We should note that for these examples, Theorem A holds. 
In fact, 
$\lambda_{\rho,K(19/45)} (1) = 1$
and
$\lambda_{\rho,K(19/45)} (-1) = 289 = 17^2$, and
$\lambda_{\rho,K(37/213)}(1)= 1$ and
$\lambda_{\rho,K(37/213)}(-1) = 225 = 15^2$.

If $p=3$, then $s_{1/3} = -1$ and for $K(r)\in H(3)$, 
$\lambda_{\rho,K(r)}(t)$
is a symmetric integer polynomial 
(of even degree) and hence, Theorem A (1)  and 
Proposition \ref{prop:3.2} (3) 
imply that  
$\lambda_{\rho,K(r)}(t)$ is the Alexander polynomial of some knot. 
Further, Theorem A (2) 
gives another condition 
that must be satisfied by this Alexander polynomial.  
Then, it is easy to show the following; 

\begin{prop}\label{prop:3.3}
The degree of $\lambda_{\rho, K(r)}(t)$ is
a multiple of $4$.
\end{prop}

\noindent
{\it Proof.} 
Write $\lambda_{\rho, K(r)}(t)=
\sum_{j=0}^{2m} a_j t^j$.
Suppose $m$ is odd, say $m=2h+~1$.
Since $\lambda_{\rho,K(r)}(1)=1$ and
$\lambda_{\rho, K(r)}(-1)=\mu^2$, it follows that
$ \sum_{j=0}^{2h} 2a_j+a_{2h+1}=1$
and $\sum_{j=0}^{2h}(-1)^j 2a_j-a_{2h+1}=\mu^2$, 
and hence, $\sum_{j=0}^{h}4a_{2j}=1+\mu^2$
that is impossible. \qed

On the other hand, for some special cases, 
it is possible to identify           
$\lambda_{\rho,K(r)}(t)$ as the Alexander polynomial of a 
certain knot.  
We can prove the following;

\begin{prop}\label{prop:3.4}
Suppose $p$ and $q$ are odd integer 
$\geq 3$.  Then 
$K(1/pq) \in H(p)$ and 
$ \lambda_{\rho,K(1/pq)}(t) =  \Delta_{K(1/q)}(t ^{2p})$,      
where  $\Delta_{K(1/q)}(t)$ is the Alexander polynomial 
of $K(1/q)$. Therefore, $\lambda_{\rho,K(1/pq)}(t)$ is
the Alexander polynomial of the $2p$-cable of the torus knot
$K(1/q)$.
\end{prop}
A proof will be given in Appendix (II).   

%

\medskip
Finally, we note that Theorem A is not true for non-rational knots.

\begin{ex}\label{ex:3.4}
Consider a non-rational knot $K =8_5$ in the Reidemeister-Rolfsen
table. 
Then $G(K)$ has a Wirtinger presentation,
$G(K) =\langle x,y,z | R_1, R_2\rangle$, 
where\\
$R_1 =  y x y^{-1}x^{-1}y^{-1}x y x y x^{-1}y^{-1}
z^{-1}y x y^{-1}x^{-1}y^{-1}z$ and\\
$
R_2 =  y x y^{-1}z y x^{-1}y^{-1} z^{-1}x^{-1}y^{-1} x z$.\\
It is easy to check that
$\rho :  x, z \mapsto
\left[\begin{array}{rr}
1&\ 1\\
0&\ 1
\end{array}\right]$ and
$y \mapsto 
\left[\begin{array}{rr}
1&\ 0\\
-1&\ 1
\end{array}\right]
$    
gives a parabolic representation of $G(K)$
on $SL(2, \ZZ)$, 
and
$\widetilde{\Delta}_{\rho,K}(t) = 
- (1-t)^2 (1+t^2) (1-2t-2t^3 -2t^5 + t^6)$. 
And hence,
       $\lambda_{\rho,K}(t) 
       =  - (1-t)^2 (1-2t-2t^3 -2t^5 + t^6)$, and 
$\lambda_{\rho,K}(1)= 0$ and  $\lambda_{\rho,K} (-1) = 2^5$.

We know,  $\lambda_{\rho,K}(t)$ is the reduced Alexander polynomial of a 
$3$-component link.
\end{ex}

\section{$\ZZ[s_0]$-Algebra}
From now on (except for Appendix), we consider exclusively the set 
$H(p), p=2n+1$.  
For simplicity, we use $s_0$ for $s_{1/p}$.   

Let $X =\xmtx$, 
and $Y =\mtx{1}{0}{s_0}{1}$
 be elements in $SL(2,\ZZ[s_0])$.\\
We define $A(x,y: \ZZ[s_0])$ as the free algebra over 
$\ZZ[s_0]$ constructed from the free group $F(x,y)$.
Let $f : A(x,y:\ZZ[s_0])  \rightarrow M_{2,2}(\ZZ[s_0])$ 
be an (algebra) homomorphism defined by  $f(x) = X$ and 
$f(y) = Y$.  
Let $S(x,y)= f^{-1}(0)$ be the kernel of $f$.   
Then $\widetilde A(s_0) = A(x,y:\ZZ[s_0])/ S(x,y)$ 
is a non-commutative $\ZZ[s_0]$-algebra.

\begin{ex}\label{ex:4.1}
 The following elements are typical elements of 
$S(x,y) :(x-1)^2$, since $(X-I)^2 = 0$, and $(y-1)^2, 
(xy)^n x (xy)^{-n}y^{-1} -1$ and $(xy)^n x-y(xy)^n$.
\end{ex}

The purpose of this section is to prove 
Lemmas \ref{lem:4.4} and \ref{lem:4.5}.
However, first we need a few technical lemmas. 

For any integer $k \geq 0$, we write

\begin{equation}
(XY)^k = \mtx{a_k}{b_k}{c_k}{d_k} 
\in SL(2,\ZZ[s_0]),
\end{equation}
where $a_k, b_k, c_k$ and $d_k$ are integer polynomials in 
$s_0$.
From the definition of $s_0$, 
we should note that $a_n=0$.



\begin{prop}\label{prop:4.1}
We have the following recursive formulas.\\
(I)\ \  \    $a_0 =  d_0 = 1$ and $b_0 = c_0 = 0$.\\
(II) \  $a_1 = 1 + s_0,  b_1 = 1,   c_1 = s_0$ and $d_1  = 1$.\\
(III)  (i)  For $k \geq 2$,
     
(1)   $a_k   = (2+s_0) a_{k-1} -  a_{k-2}$,      

(2)   $s_0 b_k = (1+s_0) a_{k-1} -  a_{k-2}$.    

\hspace*{1mm} (ii)	For $k \geq 1$,

(3)   $s_0 b_k  = a_k -  a_{k-1}$,

(4)   $s_0 b_k  = c_k$,

(5)  $a_k  = s_0 b_k + d_k$,

(6)  $d_k = a_{k-1}$,

(7)  $b_k = b_{k-1} + a_{k-1}$,

(8)  $c_k + d_k = a_k$,     

(9)   $a_0 + a_1 + \cdots + a_{k-1} = b_k$.     
\end{prop}

\noindent
{\it Proof.}
 (I) and (II) are immediate.  
 To show (III), we use induction on $k$.
For $k=1$, (3)-(9) are obvious. 
For $k = 2$, 
(1)-(9) are also 
immediate from the definition, since $(XY)^2= 
\mtx{1+3s_0+s_0^2}{2+s_0}{2s_0+s_0^2}{1+s_0}$.
Now, for any $k\geq 2$, (1) and (2) $\rightarrow$ (3), and 
(3) and (6) $\rightarrow$ (5). 
Further, since (4) and (5) $\rightarrow$ (8), and 
(2) and (3) $\rightarrow$ (7) $\rightarrow$ (9), 
it only suffices to prove (1), (2), (4) and (6).
Inductively we assume that these formulas hold for $k$.  

Then a computation $(XY)^{k+1} = (XY)^k (XY)$ shows 
\begin{align}
 &(i)\ \  \ a_{k+1} = (1+s_0) a_k + s_0 b_k, \nonumber\\ 
 &(ii)\ \ b_{k+1} =  a_k + b_k,\nonumber\\
 &(iii)\  c_{k+1} =  (1+s_0) c_{k} + s_0  d_{k},\nonumber\\
 &(iv)\  d_{k+1} = c_k + d_k.
\end{align}


And we see 

(1) $a_{k+1} = (1+s_0 ) a_k + s_0 b_k = (1+s_0)a_k + a_k
 -  a_{k-1}= (2+s_0) a_k -  a_{k-1}$,

(2) $s_0  b_{k+1} = 
s_0 a_k + s_0 b_k = s_0 a_k + (a_k - a_{k-1})
=(s_0 +1) a_k - a_{k-1}$,     

(6) $d_{k+1} = c_k + d_k = s_0 b_k + a _{k-1} =  a_k$,

(4) $c_{k+1} = (1+s_0) c_k + s_0 d_k = c_k + s_0  (c_k +  d_k)
=s_0 b_k + s_0 a_k = s_0 b_{k+1}$.

This proves Proposition \ref{prop:4.1}.
\qed

\begin{prop}\label{prop:4.2}
  (1)    $(XY)^n X=Y (XY)^n=
  \mtx{0}{b_n}{c_n}{0}$ and    
(2)      $(XY)^p = - I$.
\end{prop}

\noindent
{\it Proof.}
A direct computation shows (1), since $c_n + d_n = a_n=0$.
Also, (2) follows, since  $(XY)^p =(XY)^n X Y (XY)^n=-I$.  
\qed\\
Note that $\det\Bigl[(XY)^n X\Bigr]=- b_n  c_n=1$. 
  
\begin{prop}\label{prop:4.3}
We have the following equalities:
\begin{align}
&(1)\ a_0 + a_1 + \cdots + a_{n-1} = b_n.\nonumber\\
&(2)\ s_0 (b_1 + b_2 + \cdots + b_n) = - 1.\nonumber\\
&(3)\ b_1 + b_2 + \cdots + b_n = b_{n}^2.\nonumber\\ 
&(4)\ d_0 + d_1 + \cdots + d_n  = 1 + a_0 + a_1 + \cdots + a_{n-1}.\nonumber\\
&(5)\ c_1 + c_2 + \cdots + c_n =-1.
\end{align}
\end{prop}     

\noindent
{\it Proof.}
  First, (1) follows from Proposition 
\ref{prop:4.1}(9).
To show (2), use Proposition 
\ref{prop:4.1}(III)(3).
In fact, since $a_n (s_0) = 0$,
  $\sum_{k=1}^n s_0 b_k = \sum_{k=1}^n (a_k-a_{k-1})=a_n - a_0= -1$. 
(3) follows, since  $b_n^2 = b_n c_n s_0^{-1} = -  s_0^{-1}$. 
 (4) follows from Proposition \ref{prop:4.1} (III)(6). 
Finally, since $c_k = s_0 b_k$, it follows that
$\sum_{k=1}^n c_k = s_0   \sum_{k=1}^n  b_k = -1$, 
and (5) is proved. 
\qed

We proceed to prove two key lemmas below.
For simplicity, we use the following notations:
\begin{align}
&(1)\ 
\mbox{\boldmath $a$} = \sum_{k=0}^n a_k,  
\mbox{\boldmath $b$} = \sum_{k=0}^n b_k,  
\mbox{\boldmath $c$} = \sum_{k=0}^n c_k\  {\rm and}\  
\mbox{\boldmath $d$} = \sum_{k=0}^n d_k,\nonumber \\
&(2)\
P_k = 1 + (xy) + (xy)^2 +\cdots + (xy)^k,\nonumber\\  
&(3)\   
Q_k =  y P_k y^{-1} = 1 + (yx) + (yx)^2 + \cdots +(yx)^k.          
\end{align}

%
%
%
%

\begin{lemm}\label{lem:4.4}
 The following equalities hold in 
 $\widetilde A(s_0)$.
 \begin{align}
&(1)\ (1-y) Q_n y (1-x) = - (yx)^{n+1}(1-x).\nonumber\\
&(2)\ (1-y) Q_{2n}y (1-x) = 0.\nonumber\\
&(3)\ (1-y) Q_{3n+1}y (1-x) = - (yx)^{3n+2}(1-x). 
\end{align}
 
%
%
%
\end{lemm}

\noindent
{\it Proof.}   
Since  $(YX)^k =  Y(XY)^{k}Y^{-1}$, it follows that 
$(YX)^k = 
\mtx{d_k}{b_k}{c_k}{a_k}$.

{\it Proof of (1).}
By taking the image of both sides under $f$, we have\\
LHS $=  
\mtx{0}{0}{-s_0}{0}\mtx{\mbox{\boldmath $d$}}{\mbox{\boldmath $b$}}
{\mbox{\boldmath $c$}}{\mbox{\boldmath $a$}}
\mtx{1}{0}{s_0}{1}\mtx{0}{-1}{0}{0}$         
$=$  
$\mtx{0}{0}{0}{s_0 \mbox{\boldmath $d$}+s_0^2 \mbox{\boldmath $b$}}$.\\
On the other hand, since $x(1-x) = 1-x$ and 
$y(xy)^n = (xy)^n x$, we see \\
$ - (yx)^{n+1}(1-x) = - y(xy)^n x (1-x) = -y (xy)^n (1-x)$.  \\
Also $f(y(xy)^n) = 
\mtx{0}{b_n}{c_n}{0}$, and hence
RHS $= - 
\mtx{0}{b_n}{c_n}{0}\mtx{0}{-1}{0}{0}$
$=$
$\mtx{0}{0}{0}{c_n}$.\\                   
Use (4.3) (4) to show
$s_0 \mbox{\boldmath $d$}  +  s_0^2 \mbox{\boldmath $b$} = 
s_0 (1+ \mbox{\boldmath $a$}) -  s_0 =  s_0 \mbox{\boldmath $a$}  
= s_0 b_n = c_n$.
This proves~(1).

\medskip
{\it Proof of (2).}  Since 
\begin{align*}
(1-y) Q_{2n}y(1-x) 
&= (1-y) Q_n y (1-x) + (1-y)(Q_{2n}-Q_n) y (1-x)\\
&= - (yx)^{n+1}(1-x) + (1-y)(Q_n - 1) (yx)^n y(1-x),
\end{align*} 
it suffices to show that 
$(1-y)(Q_n - 1)(yx)^n y (1-x) =  (yx)^{n+1}(1-x)$.
Take the image of both sides under $f$. 
Then,\\
LHS $= 
\mtx{0}{0}{-s_0}{0}\mtx{\mbox{\boldmath $d$}-1}{\mbox{\boldmath $b$}}
{\mbox{\boldmath $c$}}{\mbox{\boldmath $a$}-1}\mtx{0}{b_n}{c_n}{0}
\mtx{0}{-1}{0}{0}$
$=$
$\mtx{0}{0}{0}{-c_n}$. \\
Meanwhile, RHS $= \mtx{0}{0}{0}{-c_n}$,
as is shown in the proof of (1). This proves (2).    

{\it Proof of (3).}   
Since $(yx)^{2n+1}=-1$, 
we see that
\begin{align*}
(1-y) Q_{3n+1}y(1-x) 
&= 
(1-y) Q_{2n}y(1-x) + (1-y)(Q_{3n+1} - Q_{2n})y(1-x)\\ 
&= - (1-y) Q_n y (1-x)\\
&=(yx)^{n+1}(1-x)\\
&= - (yx)^{3n+2}(1-x). 
\end{align*}
\qed

\begin{lemm}\label{lem:4.5}
The following equalities hold in $\widetilde{A}(s_0)$.
\begin{align}
 &(1)\   
(1+y) Q_n y (1+x) = (yx)^{n+1}(1+x) + 4b_n (y+ (yx)^{n+1}).
\nonumber\\
&(2)\ 
(1+y)Q_n y (1+x) (1+(xy)^n x)
= (yx)^{n+1}(1+x)(1+(xy)^n x) + 8b_n (yx)^{n+1}.
\nonumber\\
&(3)\ 
(1+y) Q_{2n}y (1+x) =  8b_n (yx)^{n+1}.
\nonumber\\
&(4)\ 
(1+y) Q_{2n}y(1+x) (1+(yx)^n y) = - 8b_n (y -  (yx)^{n+1}).
\nonumber\\
&(5)\ 
(1+y) Q_{3n+1}y (1+x) + (yx)^{n+1}(1+x)  = - 4b_n (y - (yx)^{n+1}).
\end{align}
\end{lemm}

\noindent
{\it Proof.}   
First we note that (2) follows from (1), and (4) follows 
from (3) by multiplying both sides through 
$(1+ (xy)^n x)$, since  
$(y+ (yx)^{n+1}) (1+ (xy)^n x)= y(1+ (xy)^n x)^2 = 
y (1+ 2(xy)^n x + (xy)^{2n+1}) = 2 (yx)^{n+1}$ and  
$(yx)^{n+1}(1+ (yx)^n y) = (yx)^{n+1}+ (yx)^{2n+1}y =  
(yx)^{n+1}-y$.   
Also (5) follows from (1) and (3), since 
$Q_{3n+1} = Q_{2n} - Q_n$.  
Therefore, we only need to show (1) and (3).

{\it Proof of (1).} 
Since
LHS $=(1+y)Q_n y (1+x)  = (1+y)y P_n (1+x) = y(1+y) P_n (1+x)$,
and RHS $= (yx)^{n+1}(1+x) + 4b_n y(1+ (xy)^n x) 
= y (xy)^n x (1+x) + 4b_n y(1+ (xy)^n x)$,
 it suffices to show 

(1)'  
$(1+y) P_n (1+x)  = (xy)^n x (1+x) + 4b_n (1+ (xy)^n x)$.

By taking the image of both sides of (1)' under $f$,
we obtain

\begin{align*}
{\rm LHS} &= \mtx{2}{0}{s_0}{2}\mtx{\mbox{\boldmath $a$}}
{\mbox{\boldmath $b$}}{\mbox{\boldmath $c$}}{\mbox{\boldmath $d$}}\mtx{2}{1}{0}{2}=       
\mtx{4\mbox{\boldmath $a$}}{2\mbox{\boldmath $a$}+4\mbox{\boldmath $b$}}
{2s_0 \mbox{\boldmath $a$}+4\mbox{\boldmath $c$}}{s_0 \mbox{\boldmath $a$}
+2\mbox{\boldmath $c$}+2s_0 \mbox{\boldmath $b$} +4\mbox{\boldmath $d$}},\ {\rm and}\\
 {\rm RHS}  &= 
\mtx{0}{b_n}{c_n}{0}\mtx{2}{1}{0}{2}+
4b_n\mtx{1}{b_n}{c_n}{1}=\mtx{4b_n}{2b_n+4b_n^2}{2c_n+4b_n c_n}{c_n+4b_n}.
\end{align*}

Therefore we need to show
\begin{align}
&(i)\
4 \mbox{\boldmath $a$} = 4 b_n,
\nonumber\\ 
&(ii)\
2 \mbox{\boldmath $a$}  + 4 \mbox{\boldmath $b$}  
= 2b_n + 4b_n^2,
\nonumber\\ 
&(iii)\
2 s_0 \mbox{\boldmath $a$} + 4 \mbox{\boldmath $c$}  
=  2c_n + 4 b_n c_n,\ {\rm and}
\nonumber\\
&(iv)\ 
s_0 \mbox{\boldmath $a$} + 2 \mbox{\boldmath $c$} 
+ 2 \mbox{\boldmath $b$} s_0 + 4 \mbox{\boldmath $d$} = c_n + 4b_n.
\end{align}
%
%
%
%

First, (i) follows from (4.3)(1), 
and (ii) follows from (4.3)(1) and (3). 
Further, (iii) follows, 
since  $s_0 \mbox{\boldmath $a$}  = s_0 b_n  = c_n$  and  
$\mbox{\boldmath $c$}  = s_0 \mbox{\boldmath $b$}  = -1 =b_n c_n$.
Finally, (iv) follows, since 
$s_0 \mbox{\boldmath $a$} + 2 \mbox{\boldmath $c$} 
= c_n-2$, and
$2 \mbox{\boldmath $b$} s_0 + 4 \mbox{\boldmath $d$} 
= -2 + 4(1 + \mbox{\boldmath $a$} ) = 2 + 4b_n$
by (4.3) (1) and~(5).
A proof of (1) is now complete.

{\it Proof of (3).}  First, we note\\
LHS $= (1+y) Q_{2n}y(1+x)  = (1+y)Q_n y(1+x) 
 + (1+y)(Q_{2n}-Q_n)y(1+x)$.\\
Since by (4.6)(1), 
$(1+y)Q_n y(1+x)= (yx)^{n+1}(1+x) + 4b_n (y+(yx)^{n+1})$ and \\
$(1+y)(Q_{2n}-Q_n)y(1+x) = (1+y)(Q_n - 1) (yx)^n y(1+x)$,
we must show\\
$(yx)^{n+1}(1+x) + 4b_n (y+(yx)^{n+1}) + (1+y)(Q_n -1)(yx)^n y(1+x) 
= 8b_n (yx)^{n+1}$.\\
This equation is equivalent to  \\
 $(yx)^{n+1}(1+x) + (1+y)(Q_n - 1) (yx)^n y(1+x) 
 = 4b_n y\{(xy)^n x -1\}$. \\
But since $Q_n = y P_n y^{-1}$, it suffices to show
\begin{equation}
\noindent
(xy)^n x (1+x) + (1+y)(P_n - 1) (xy)^n (1+x) 
 = 4b_n \{(xy)^n x -1\}.
\end{equation}


Now take the image of both sides of (4.8) under $f$.  
Then,
\begin{align*}
{\rm LHS} &= 
\mtx{0}{b_n}{c_n}{0}\mtx{2}{1}{0}{2}
+\mtx{2}{0}{s_0}{2}\mtx{\mbox{\boldmath $a$}-1}{\mbox{\boldmath $b$}}
{\mbox{\boldmath $c$}}{\mbox{\boldmath $d$}-1}\mtx{0}{2b_n}{2c_n}{d_n}\\
&=
\left[\begin{array}{ll}
4\mbox{\boldmath $b$}c_n
&\ \  \ 4\mbox{\boldmath $a$}b_n-2b_n+
2\mbox{\boldmath $b$}d_n\\
4c_n(\mbox{\boldmath $d$}-1)
&\ \ \  c_n+2b_n s_0(\mbox{\boldmath $a$}-1)+4b_n \mbox{\boldmath $c$}+
d_n(2\mbox{\boldmath $d$}-3)
\end{array}
\right],\ {\rm and}\\
{\rm RHS}  &=  4b_n 
\mtx{-1}{b_n}{c_n}{-1}.
\end{align*}

\noindent      
Therefore, we need to show 
\begin{align}
&(i)\    
4 \mbox{\boldmath $b$}c_n =  - 4 b_n,
\nonumber\\ 
&(ii)\
4 \mbox{\boldmath $a$} b_n -2b_n + 2\mbox{\boldmath $b$}
d_n =  4b_n^2,
\nonumber\\
&(iii)\
4c_n (\mbox{\boldmath $d$}-1)  =  4 b_n c_n,
\nonumber\\
&(iv)\
c_n + 2s_0 b_n (\mbox{\boldmath $a$} -1) + 4b_n 
\mbox{\boldmath $c$}  + d_n (2\mbox{\boldmath $d$} -3) = -4b_n.
\end{align}      

First, (i) follows, since  $4\mbox{\boldmath $b$}c_n = 4\mbox{\boldmath $b$}s_0 b_n = - 4 b_n$, 
and (ii) follows, 
since $4\mbox{\boldmath $a$}b_n - 2b_n +2\mbox{\boldmath $b$}d_n 
= 4 b_n^2 - 2b_n - 2\mbox{\boldmath $b$}s_0 b_n  
= 4b_n^2 -2b_n + 2b_n= 4b_n^2$.   
Note that $d_n = a_{n-1} = - s_0 b_n$.   
 Also, (iii) follows, since  $\mbox{\boldmath $d$}-1=b_n$.   
 Finally, (iv) follows, since  \\
$c_n +2s_0 b_n (\mbox{\boldmath $a$}-1) + 4b_n \mbox{\boldmath $c$} + d_n (2\mbox{\boldmath $d$}-3) = 
s_0 b_n +2s_0 b_n (b_n - 1)  - 4b_n  + 
+ a_{n-1}(2\mbox{\boldmath $a$}-1) = s_0 b_n +2s_0 b_n^2 - 2s_0 b_n -4b_n-
s_0 b_n (2b_n - 1) = -4b_n$.  \\
This proves (3). \qed
%
%
%
%
%
%
%
%
%
%

\section{Restatement of Theorem A.}

Let $K(r)$ be an element of $H(p)$. 
Then  $r = \frac{\beta}{\alpha}$ has a continued fraction expansion of the form:
 $\frac{\beta}{\alpha} = 
 [pk_1,2m_1, pk_2, 2m_2,\ldots,  2m_q, pk_{q+1}]$,
where  $k_i, m_j$ are non-zero integers.  \\
Using this form, we can construct a diagram of $K(r)$ as a 
$4$-plat.
First construct a $3$-braid 
$\gamma =\sigma_2^{pk_1}\sigma_1^{2m_1}\sigma_2^{ pk_2}\cdots \sigma_1^{2m_q}\sigma_2^{pk_{q+1}}$,
where $\sigma_i$ are Artin's generators of the $3$-braid group. 
See Fig 5.1.  
    
\begin{figure}[h]
\begin{center}
\includegraphics[height=12mm]{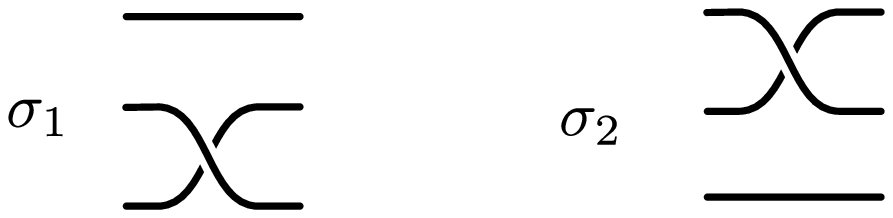}
\end{center}
\centerline{Figure 5.1:
The Artin generators for $3$-braids}
\end{figure}

\begin{figure}[h]
\begin{center}
\includegraphics[height=15mm]{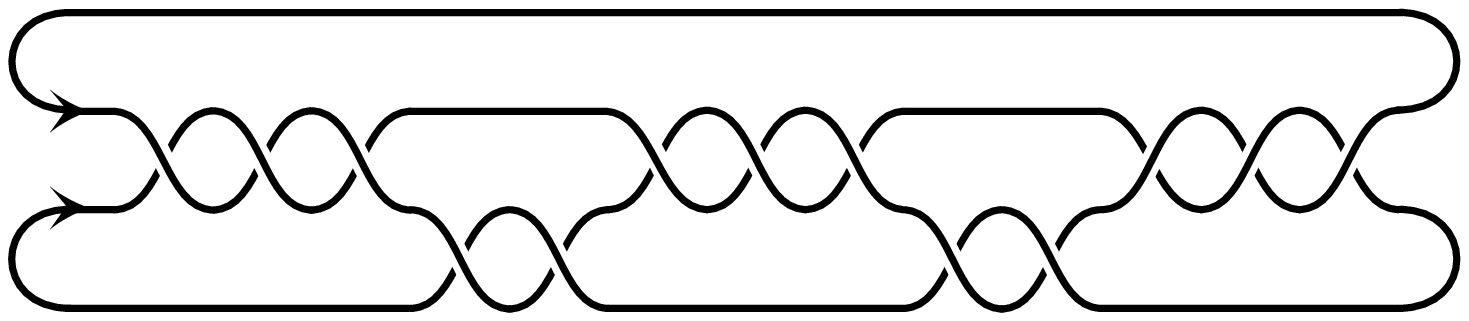}
\end{center}
\centerline{Figure 5.2: 
$2$-bridge knot of $29/69 = [3,2,3,2,-3]$}
\end{figure}

Close $\gamma$ by joining the first and second strings 
(at the both ends) and then join the top and 
bottom of the third string by a simple arc as 
in Fig. 5.2.   (For convenience, figures will be $\pi/2$ rotated.)  
We give downward orientation to
the second and third strings.

Fig 5.2 shows the 
(oriented) $2$-bridge knot obtained from the 
continued fraction  
$29/69 = [3,2,3,2,-3]$.  
(A braid gives a knot diagram $D(r)$ of $K(r)$ if and only if 
$\sum_{j=1}^{q+1}k_j$ is odd.)

\begin{rem}\label{rem:5.1}
Although $k_i$ and $m_j$ are not $0$, later in this paper,
we need an appropriate interpretation of our continued
fractions when some are $0$. 
The following interpretations will be easily justified 
by checking their diagrams as $4$-plats.
If $k_i\ (i\neq 1, q+1)$ or $m_j\ (1\le j\le q)$ is $0$,
then $r$ is interpreted as \\
$[pk_1,2m_1,\dots, pk_{i-1}, 2(m_{i-1}+m_i),pk_{i+1},\dots,pk_{q+1}]$
or\\
$[pk_1,2m_1,\dots,
2m_{j-1},p(k_j+k_{j+1}),2m_{j+1},\dots, pk_{q+1}]$.
If $k_1=0$ (or $k_{q+1}=0$),
then $r$ is interpreted as 
$[pk_2,2m_2,\dots, pk_{q+1}]$
(or $[pk_1,2m_1,\dots, pk_q]$).
Note that our continued fraction expansions start and end with $pk$.

\end{rem}

Next we find a presentation of $G(K)$ from $D(r)$.  
Two (meridian) generators $x$ and $y$ are represented 
by loops that go around once under local maximal points 
from the left to the right as shown in Fig.5.3.  
The relation is obtained using $x$ and $y$ by a standard method.
However, we describe this process more precisely.

\begin{figure}[h]
\begin{center}
\includegraphics[height=30mm]{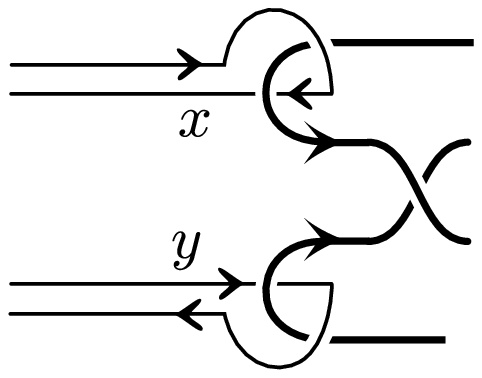}
\end{center}
\centerline{Figure 5.3: Generators for the knot group}
\end{figure}

First we divide $D(r)$ into $2q+3$ small pieces by 
$2q+2$ vertical lines $L_j$.  
See Fig 5.4. 
We define the elements $x_0, x_1,\cdots, x_{2q+1}, y_0, y_1, 
\cdots, y_{2q+1}$, $z_0, z_1, \cdots, z_{2q+1}$ in  
$F(x,y)$ as follows.   
Let $Z_j, X_j, Y_j$, respectively, be the points of 
intersection of $L_j$ and the first, second and third strings.  
Then  $z_j, x_j, y_j$, respectively, are represented by loops 
that go around once under these points $Z_j, X_j, Y_j$ from the 
left to the right.  
We note that $x_0 = x ,  y_0 = y$ and $z_0 = x$.  

\begin{figure}[h]
\includegraphics[height=37mm]{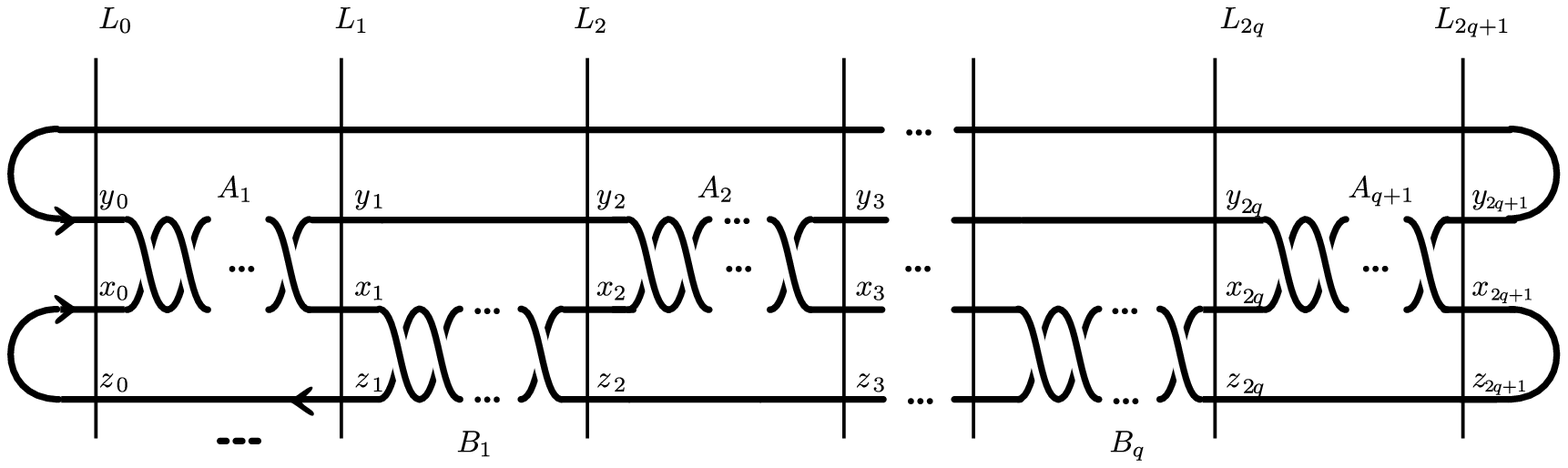}
\centerline{Figure 5.4: Elements in $F(x, y)$}
\end{figure}

By the standard method, $x_j, y_j, z_j$ can be written as 
words of $x$ and $y$. See Fig. 5.5.
Let $A_{j+1} = y_{2j}x_{2j}$, $j=0, 1, \cdots, q$. 
Then  $x _{2j+1}$ and $y_{2j+1}$ are given as follows:

\begin{align}
&(1)\  {\rm If}\ 
k_j = 2\ell_j,\ {\rm then}\
y_{2j+1} = A_{j+1}^{p\ell_j} y_{2j}A_{j+1}^{-p\ell_j},
\ {\rm and}\  
x_{2j+1} = A_{j+1}^{p\ell_j} x_{2j}A_{j+1}^{-p\ell_j}.
\nonumber\\
&(2)\ {\rm If}\  
k_j = 2\ell_j + 1,\ {\rm then}\ 
y_{2j+1} = A_{j+1}^{p\ell_j+n+1}
 x_{2j}A_{j+1}^{-(p\ell_j+n+1)},{\rm and}\  
 \nonumber\\
 &\hspace{3.9cm}
x_{2j+1} 
= A_{j+1}^{p\ell_j+n} y_{2j}
A_{j+1}^{-(p\ell_j+n)}.
\nonumber\\
&(3)\ 
y_{2j+2} = y_{2j+1}\ {\rm and}\   
z_{2j+1} = z_{2j}, j = 0, 1,
\cdots,q.
\end{align}

%
%
%
%

Let  $B_{j+1} = x_{2j+1} z_{2j+1}^{-1}$, $j=0, 1, \cdots,q-1.$ 
Then  $x_{2j+2}$ and $z_{2j+2}$ are given by
   
\begin{equation}
x_{2j+2} = B_{j+1}^{m_j}x_{2j+1}B_{j+1}^{-m_j},\  {\rm and}\  
z_{2j+2} = B_{j+1}^{m_j}z_{2j+1}B_{j+1}^{-m_j}.
\end{equation} 

\begin{figure}[h]
\begin{center}
\includegraphics[height=40mm]{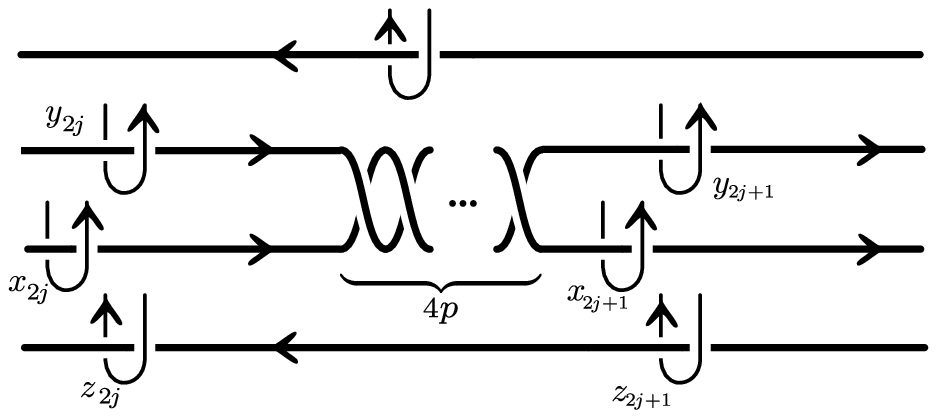}
\end{center}
\centerline{Figure 5.5: Rewriting process of letters}
\end{figure}

Then the relation of $G(K(r))$ is given by

%
\begin{equation}
y_{2q+1} =  y,\ {\rm  (or equivalently}\   
x_{2q+1} = z_{2q+1}.)
\end{equation}

Therefore, $G(K(r)) = \langle x,y|R\rangle$
 is a Wirtinger presentation of 
$G(K(r))$, where $R = y_{2q+1}y^{-1}$.
We note that relation (5.3) is a conjugate of the 
relation given by (2.2).

Now we can express the relation $y_{2q+1}y^{-1}$ as a 
product of conjugate of  $R_0 = (xy)^n x (xy)^{-n} y^{-1}$:

\begin{equation}
y_{2q+1}y^{-1}= \prod_{j=1}^m u_j 
R_0^{\epsilon_j} u_j^{-1}, \ 
{\rm where}\  u_j \in F(x,y)\ {\rm  and}\  
\epsilon_j = \pm 1.
\end{equation}

Let  $\Phi_0:\widetilde A(s_0) \rightarrow M_{2,2}(\ZZ[s_0])$ 
be a homomorphism defined by 
$\Phi_0 = \Phi|_{t=1}$, and hence $\Phi_0 (x) = \rho (x)$ and
$\Phi_0 (y) = \rho (y)$.  Then it follows from (3.4) that
%

\begin{align}
&(1)\    
\lambda_{\rho,K(r)}(1)=
\det \left[\sum_{j=1}^m \epsilon_j u_j^{\Phi_0}\right],\ {\rm and}\ 
\nonumber\\ 
&(2)\  
\lambda_{\rho,K(r)}(-1)=
\det \left[\sum_{j=1}^m (-1)^{\ell(u_j)}\epsilon_j u_j^{\Phi_0}\right],
\end{align}
where $\ell(u_j)$ denotes the length of a word $u_j \in F(x,y)$. \\
Therefore, to prove Theorem A, it will be sufficient to show the following proposition.

\begin{prop}\label{prop:5.1}
(1) The element  
$\lambda (r)= \sum_{j=1}^m \epsilon_j u_j$ in the 
$\ZZ[s_0]$-algebra $\widetilde A(s_0)$ is a single element, 
namely, 
 $\lambda (r)=\sum_{j=1}^m \epsilon_j u_j=\pm w$,
for some element $w$ in $F(x,y)$.

(2)   The element  $\widetilde\lambda (r)= \sum_{j=1}^m 
(-1)^{\ell(u_j)}\epsilon_j u_j$  is a constant multiple of 
a single element, i.e.
$\widetilde\lambda (r) =\pm\mu w$ for some  
$\mu \in \ZZ[s_0]$ and $w \in F(x,y)$.
\end{prop}

\section{Rewriting process}

Now $R_0=(xy)^n x(xy)^{-n}y^{-1}$ 
is a defining relation of $G(K(1/p))$. 
We denote $uR_0^{\epsilon}u^{-1}$ by 
$R_0^{\epsilon u}$, for $u \in F(x,y)$ 
and $\epsilon =\pm 1$. 

In this section, we establish a rewriting process which transforms 
an element $w \in F(x,y)$ into the form $R_0^uw_0$, 
where $u \in \widetilde A(s_0)$ and $w_0 \in F(x,y)$.  
Since we are concerned on an element $\lambda (r)$ or 
$\widetilde\lambda (r)$ of $\widetilde A(s_0)$,
we may write  $R_0^u = R_0^{u'}$ if $u = u'$
in  $\widetilde A(s_0)$, 
and $R_0^u R_0^v = R_0^{u+v}= R_0^v R_0^u$, where 
$u, v \in \widetilde A(s_0)$.



\begin{lemm}\label{lem:6.1}
We have the following formulas involving $R_0$.\\
(I) 
(1) $(yx)^{n+1}x (yx)^{-(n+1)}= R_0^y y$.\\
(2) $(yx)^n y (yx)^{-n}= R_0^{-1}x$.\\
(3) $(yx)^{2n+1}x (yx)^{-(2n+1)}= 
R_0^{(yx)^n y-1}x = R_0^{(xy)^n x-1}x$.\\
(4) $(yx)^{2n+1}y (yx)^{-(2n+1)}= 
R_0^{-(yx)^{n+1}+ y}y$.\\
(5) $(yx)^{3n+2}x (yx)^{-(3n+2)} =
R_0^{(yx)^{p}y-(yx)^{n+1}+y}y= R_0^{-(yx)^{n+1}}y$. \\   
(6) $(yx)^{3n+1}y (yx)^{-(3n+1)} = 
R_0^{-(yx)^{p}+(yx)^{n}y-1}x
=
R_0^{(yx)^n y}x=
R_0^{x(yx)^n}x$.\\
(II)  For $k \geq 1$, 
$(R_0^g u)^k = R_0^{(1+u+\cdots+u^{k-1})g} u$, where  
$g \in \widetilde A(s_0)$ and $u\in F(x,y)$.
\end{lemm}

\noindent
{\it Proof.} Since  most of our proofs are straightforward, 
we prove only one of these formulas, say (I) (3).
In fact, since $(yx)^n y = (xy)^n x$ in $\widetilde A(s_0)$,
we have:
\begin{align*}
(yx)^{2n+1}x (yx)^{-(2n+1)} 
&=(yx)^{2n+1}(y^{-1}   x^{-1})^{2n}y^{-1} \\
&=(yx)^n y(xy)^n x (y^{-1}x^{-1})^n y^{-1} 
(x^{-1}y^{-1})^n \\
&=(yx)^n y R_0 (x^{-1}y^{-1})^n\\
&=R_0^{(yx)^n y}(xy)^n y (x^{-1}y^{-1})^n x^{-1} x\\
&=R_0^{(yx)^n y} R_0^{-1} x  \\
&= R_0^{(yx)^n y -1}x\\ 
&=  R_0^{(xy)^n x-1}x.
\end{align*} 
 \qed

\begin{lemm}\label{lem:6.2}
For the elements defined in Section 5,
we have\\
(1) 
For  $j = 0,1,2,\cdots, 2q +1,
y_j x_j z_j = y$, as elements of  $F(x,y)$.\\
(2) 
For  $j = 0,1,2,\cdots, 2q+1$,
we can write 
$x_j = R_0^{v_j}x,  y_j = R_0^{w_j}y$ and  
$z_j = R_0^{u_j}x$, where 
$u_j, v_j, w_j$ are elements of  
$\widetilde A(s_0)$.

%
\end{lemm}

\noindent
{\it Proof.}  
(1) is evident from the definition of $x_j, y_j, z_j$.
Further for $j = 0$, (2) is evident.  
Consider the case $j=1$.

If $k_1 =2\ell_1$, 
we apply Lemma \ref{lem:6.1}
(I)(4) repeatedly to obtain \\
$y_1 = (yx)^{p\ell_1}y (yx)^{-p\ell_1} =  R^w y$   
for some $w \in \widetilde A(s_0)$.
If $k_1=2\ell_1 +1$, 
then by Lemma \ref{lem:6.1}
(I)(3) we see that
\begin{align*}
y_1 &= (yx)^{p\ell_1+n+1}x (yx)^{-(p\ell_1+n+1)}\\
&= (yx)^{n+1}(yx)^{p\ell_1}x (yx)^{-p\ell_1}(yx)^{(n+1)}\\
&= (yx)^{n+1}R_0^v x (yx)^{(n+1)}
\ \ ({\rm for\ some}\ v \in \widetilde A(s_0))\\
&=  R_0^{(yx)^{n+1}v} (yx)^{n+1}x (yx)^{-(n+1)} \\
&=  R_0^{(yx)^{n+1} v+y} y  \ \ ({\rm by\ Lemma}\
 \ref{lem:6.1}{\rm (I)}(1)).)
 \end{align*}
Also, $z_1  =  z_0  = x = R_0^0 x$.
Using Lemma \ref{lem:6.1}
(II), we can complete the proof by an easy inductive argument.  
The details are omitted.
\qed

Note that to prove Theorem A we need more precise 
description of these elements $w_j, u_j, v_j$ 
that will be given in 
Proposition \ref{prop:7.1}.

Now Lemma \ref{lem:6.2} 
makes our proof of Theorem A considerably simpler
as shown in the following proposition.

\begin{prop}\label{prop:6.3}
Let $r = 
[pk_1,2m_1, pk_2, 2m_2, \cdots, 2m_q, pk_{q+1}]$ and \\
$r' = 
[pk'_1,2m_1, pk'_2, 2m_2, \cdots, 
2m_q, pk'_{q+1}]$.
Then $\lambda (r) =  
\lambda (r')$  if  $k_j  \equiv  k'$   (mod 4)\ 
for\ $j = 1,2, \cdots, q +1$.
\end{prop}

We should note that even though 
$\lambda (r) =  \lambda (r')$, 
their twisted Alexander polynomials are different.

\medskip
\noindent
{\it Proof of Proposition \ref{prop:6.3}.} 
Suppose that $k_j = 4$. 
By Lemma \ref{lem:6.2}, 
we can write $y_{2j} = y_{2j-1} = R_0^w y$ and  
$x_{2j} = R_0^v x$ for some $w,v \in \widetilde A(s_0)$.   
Since  $A_{j+1} = y_{2j} x_{2j}=R_0^w y R_0^v x 
= R_0^{w+yv}(yx)$, it follows that 
\begin{align*}
y_{2j+1}&= A_{j+1}^{2p}y_{2j}A_{j+1}^{-2p} \\
&= (R_0^{w+yv}yx)^{2p}y_{2j} (R_0^{w+yv}yx)^{-2p}.
\end{align*}

Let $w+yv=g$.  Then, by Lemma \ref{lem:6.1}
(II), we see \\
\begin{align*}
y_{2j+1} &= (R_0^g yx)^{2p} y_{2j}(x^{-1}y^{-1}R_0^{-g})\\
& = R_0^{Q_{2p-1}g} (yx)^{2p}R_0^w y (x^{-1} y^{-1})^{2p}
 R_0^{-Q_{2p-1} g}.
 \end{align*}
However, $Q_{2p-1}=0$, since $(yx)^p=-1$, by Proposition
\ref{prop:4.2} (2),
and hence,     
\begin{align*}
y_{2j+1} &= (yx)^{2p}R_0^w y(x^{-1}y^{-1})^{2p} \\
&= R_0^  {(yx)^{2p}w} (yx)^{2p} y (x^{-1}y^{-1})^{2p}\\
&= R_0^ w (yx)^{2p} y (x^{-1}y^{-1})^{2p}.
\end{align*}
By using Lemma \ref{lem:6.1}
(I)(4), we can show 
$(yx)^{2p} y (x^{-1}y^{-1})^{2p} 
= R_0^{(yx)^p y +y} y=  y$,
and hence, we have $y_{2j+1} = R_0^w y$.
Similarly, we obtain $x_{2j+1} = R_0^v x$.
The same argument works for $k_j=-4$.\\
Now we know that we may replace 
$\sigma_2^{pk_j}$ by  $\sigma_2^{pk_j \pm 4p}$  
in the braid presentation $\gamma$  of $K(r)$ defined 
in Section 5 keeping 
$\lambda (r)$  unchanged.  
This proves Proposition \ref{prop:6.3}.
\qed

By Proposition \ref{prop:6.3}, 
we may assume that
\begin{equation}
k_j = 1,2\ {\rm  or}\ 3\  {\rm for\  any\ } 
j=1,2,\cdots, q +1.
\end{equation}
If $k_j \equiv 0$ (mod 4), then we may take $k_j = 0$ and  
$r$ is reduced to a shorter continued fraction
(Remark \ref{rem:5.1}).

\section{Proof of Theorem A (I),  Proof of Proposition 5.2 (1)}

Let  $r  =[pk_1,2m_1, pk_2, 2m_2, \cdots, 2m_q, pk_{q+1}]$,
and
we may assume that $k_j = 1, 2$ or $3$ for 
$1 \leq j \leq q+1$ and $m_j \ne 0$  for $1 \leq j \leq q$.\\
We want to show that 
\begin{equation}
\lambda (r)=\sum_{j=1}^m \epsilon_j u_j=\pm w,\    
{\rm for\  some\ } w \in F(x,y).
\end{equation}

First we determine precisely the elements $w_j, u_j, v_j$.

\begin{prop}\label{prop:7.1}
For any $j =  0, 1, 2, \cdots,2q +1$,
we can write
$y_j = R_0^{w_j}y, x_j = R_0^{v_j}x$  and  
$z_j = R_0^{u_j} x$.
Then we have the following;\\
(1)   $u_0 = u_1 = 0,  v_0 = 0,  w_0 = 0$.\\
 (2)   $u_0 = u_1, u_2 = u_3, \cdots, u_{2q} = u_{2q+1}$. \\     
(3)   $w_1= w_2, w_3 = w_4, \cdots, w_{2q-1}= w_{2q}$.  \\
(4)    $u_{2j} = u_{2j+1} = \sum_{k=1}^j m_k (x-1) y^{-1} 
w_{2k-1}, j = 1,2, \cdots, q$.\\
(5)     $v_j = u_j - y^{-1}w_j, j=1,2, \cdots,2q +1$.\\
(6)   For $j= 0,1,2, \cdots, q$, 
%
%
%

$$w_{2j+1}=
\begin{cases}
0&{\rm if\ }\sum_{i=1}^{j+1} k_i  \equiv 0 \ ({\rm mod\ }4) \\
y&{\rm if\ }\sum_{i=1}^{j+1} k_i  \equiv 1 \ ({\rm mod\ }4) \\
 y-(yx)^{n+1}&
 {\rm if\ }\sum_{i=1}^{j+1} k_i  \equiv 2\ ({\rm mod\ }4) \\
- (yx)^{n+1} &
{\rm if\ } \sum_{i=1}^{j+1} k_i  \equiv 3
 \ ({\rm mod\ }4) \\
\end{cases}
$$

\end{prop}

Since  
$\lambda (r) = w_{2q+1}$, 
it follows that if $K(r)$ is a knot, 
then $ \lambda (r) = y$ or  $-(yx)^{n+1}$. 
 This proves (7.1) and hence Proposition \ref{prop:5.1} (1). 

\noindent
{\it Proof of Proposition \ref{prop:7.1}.}
Formulas (1)-(3) follow from the diagram $D(r)$.  
Further, (5) follows 
from Lemma \ref{lem:6.2}(1).
In fact, $y^{-1}y_j x_j z_j^{-1}= 1$
yields\\
$y^{-1}R_0^{w_j} y R_0^{v_j} x 
x^{-1}R_0^{-u_j}= 
R_0^{y^{-1}w_j}R_0^{v_j-u_j} = 1$ 
and hence
$y^{-1}w_j + v_j - u_j =0$ and  $v_j = u_j - y^{-1}w_j$.

Now we prove (4) by induction on $j$.
Since $B_1 = x_{1}z_{1}^{-1} = R_0^{v_1}x x^{-1} = 
R_0^{v_1}$, we see that 
\begin{equation}
z_2 = B_ 1^{m_1} z_1B_1^{-m_1} = 
R_0^{m_1 v_1} x  R_0^{-m_1v_1} = R_0^{m_1(1-x)v_1} x.
\end{equation}  

Therefore, $u_2 = m_1(1-x)v_1$.  
Since $v_1 = u_1 -y^{-1}w_1 = - y^{-1}w_1$, 
we have   $u_2 = m_1(x-1) y^{-1}w_1$.   
This proves $(4)_1$.

Next consider $u_{2\ell+2}$.
By induction, we assume that \\
$u_{2\ell} = \sum_{j=1}^\ell m_j (x-1) y^{-1}w_{2j-1}$. \\
Since $B_{\ell+1} = x_{2\ell+1}z_{2\ell+1}^{-1} =  
R_0^{v_{2\ell +1}}x x^{-1}R_0^{-u_{2\ell}} = 
R_0^{v_{2\ell +1} -u_{2\ell}}                   
= R_0^{-y^{-1}w_{2\ell +1}}$,       
we have 
\begin{align}
z_{2\ell+2} &= B_{\ell+1}^{m_{\ell+1}} 
z_{2\ell}B_{\ell+1}^{-m_{\ell+1}}   \nonumber\\      
 &= R_0^{-m_{\ell+1}y^{-1}w_{2\ell+1}}R_0^{u_{2\ell}} x 
R_0^{m_{\ell+1}y^{-1}w_{2\ell+1}}\nonumber\\
&= R_0^{-m_{\ell+1}(1-x)y^{-1}w_{2\ell+1}+u_{2\ell}}x.
\end{align}


Since $u_{2\ell+1} = u_{2\ell}$, we see
\begin{align*}
u_{2\ell+2} 
&= - m_{\ell+1}(1-x) 
y^{-1}w_{2\ell+1} + \sum_{j=1}^{\ell}m_j 
(x-1) y^{-1}w_{2j-1}  \\ 
&= \sum_{j=1}^{\ell +1}m_j (x-1) y^{-1}w_{2j-1}.
\end{align*}   
This proves $(4)_{\ell+1}$.  

Finally we prove (6) by induction on $j$. 
Consider the initial case $w_1$.

Case 1.  $k_1 = 1$.
Since $A_1 = y_0 x_0 = yx$, we see from Lemma \ref{lem:6.1}(I)(1),
\begin{equation}
y_1 = A_1^{n+1}x_0 A_1^{-(n+1)} = 
(yx)^{n+1}x (x^{-1}y^{-1})^{n+1} = R_0^y y. 
\end{equation}
Therefore,  $w_1 = y$.  

Case 2.   $k_1 = 2$.   
As is seen in Case 1, we have from Lemma \ref{lem:6.1}(I)(4)
\begin{equation}
y_1 = A_1^p y_0 A_1^{-p} = R_0^{-(yx)^{n+1}+y}y, 
{\rm \ and\  hence}\ w_1 = -(yx)^{n+1}+ y.
\end{equation} 

Case 3.   $k_1 = 3$.   
Then $pk_1 = 6n+3$ and 
\begin{equation}
y_1 = A_1^{3n+2}x_0 A_1^{-(3n+2)} = 
(yx)^{3n+2}x (x^{-1}y^{-1}) ^{3n+2}        
= R_0^{-(yx)^{n+1}} y, 
\end{equation}
and hence, $w_1 =  -(yx)^{n+1}$.  This proves $(6)_1$.    

Next consider  $w_{2\ell+1}$.  
Again the proof is divided into three cases: 
$k_{\ell+1}= 1,2,3$.

Case 1.   $k_{\ell+1} = 1$ and $pk_{\ell+1} = 2n+1$.\\
Then, $A_{\ell+1} = y_{2\ell}x_{2\ell} = R_0^{w_{2\ell}} 
y R_0^{v_{2\ell}} x = R_0^{w_{2\ell}+yv_{2\ell}}yx$,  and 
 $y_{2\ell+1} = A_{\ell+1}^{n+1} x_{2\ell}A_{\ell+1}^{-(n+1)}$.

By induction assumption, we have:\\
$w_{2\ell} + y v_{2\ell} 
= w_{2\ell} + y(u_{2\ell}-y^{-1}w_{2\ell}) 
= y u_{2\ell}$.
Therefore,  
$A_{\ell+1} = R_0^{y u_{2\ell}} yx$  and hence
$y_{2\ell+1} = (R_0^{y u_{2\ell}} yx)^{n+1} 
R_0^{v_{2\ell}}x (R_0^{y u_{2\ell}} yx)^{-(n+1)}$. \\
Since  $A_{\ell+1}^{n+1} = 
R_0^{Q_n y u_{2\ell}}(yx)^{n+1}$, it follows that
\begin{align*}
y_{2\ell+1}  &= R_0^{Q_n y u_{2\ell}}(yx)^{n+1}
R_0^{v_{2\ell}}x (yx)^{-(n+1)} R_0^{-Q_n y u_{2\ell}}\\
  & =R_0^{Q_n y u_{2\ell}} R_0^{(yx)^{n+1}v_{2\ell}}
(yx)^{n+1}
x(yx)^{-(n+1)} R_0^{-Q_n y u_{2\ell}}\\
  & = R_0^{Q_n y u_{2\ell}+ (yx)^{n+1}v_{2\ell}+y} y  
R_0^{-Q_n y u_{2\ell}}\\
  & = R_0^{(1-y)Q_n y u_{2\ell}+ (yx)^{n+1}v_{2\ell}+y} y 
  \end{align*}
                                                      
and hence
\begin{equation}
w_{2\ell+1} 
= (1-y)Q_n yu_{2\ell}+ (yx)^{n+1}v_{2\ell}+y.
\end{equation}

Since by (4), 
$u_{2\ell} = \sum_{j=1}^{\ell} m_j (x-1) 
 y^{-1} w_{2j-1}$ and 
 $v_{2\ell} = u_{2\ell}-y^{-1}w_{2\ell-1}$, we have 
$w_{2\ell+1} = $
$
\{(1-y)Q_n y +(yx)^{n+1}\}
\left(
{\displaystyle
\sum_{j=1}^{\ell}} m_j (x-1) y^{-1} w_{2j-1}\right) 
- (yx)^{n+1} y^{-1}w_{2\ell-1}+ y$. \\
But $\{(1-y) Q_n y + (yx)^{n+1}\}(x-1) = 0$, 
by (4.5)(1), and hence, 
\begin{equation}
w_{2\ell+1} 
=  -(yx)^{n+1}y^{-1}w_{2\ell-1} + y.
\end{equation}

Now, we consider the following four subcases separately.\\
Case (i)  $\sum_{j=1}^{\ell} k_j \equiv 0$  (mod 4), and thus 
$\sum_{j=1}^{\ell+1} k_j \equiv 1$ (mod 4).\\
Then by induction $w_{2\ell-1}= 0$ and hence 
$w_{2\ell+1}= y$.     \\
Case (ii) $\sum_{j=1}^{\ell} k_j \equiv 1$  (mod 4), 
and thus $\sum_{j=1}^{\ell+1} k_j \equiv 2$ (mod 4).\\
Then by induction $w_{2\ell-1} = y$ and hence 
$w_{2\ell+1} = -(yx)^{n+1} + y$. \\
Case (iii)   $\sum_{j=1}^{\ell} k_j \equiv 2$  (mod 4), 
and $\sum_{j=1}^{\ell+1} k_j \equiv 3$ (mod 4).\\
Since $w_{2\ell-1} = y - (yx)^{n+1}$,
we have 
\begin{align*}
w_{2\ell+1} 
&=  -(yx)^{n+1}      y^{-1}(y-(yx)^{n+1}) + y\\
&=  -(yx)^{n+1} + (yx)^{n+1}(yx)^n y + y \\
&=  -(yx)^{n+1}.   
\end{align*}
Case (iv)   $\sum_{j=1}^{\ell} k_j \equiv 3$  (mod 4), 
and $\sum_{j=1}^{\ell+1} k_j \equiv 0$ (mod 4).\\
Since $w_{2\ell-1} = - (yx)^{n+1}$, 
we have 
\begin{align*}
w_{2\ell+1} 
&=  -(yx)^{n+1} y^{-1}(-(yx)^{n+1})~+~y\\
&= (yx)^{n+1}(xy)^n x + y = (yx)^{n+1}(yx)^n y + y \\
&=  -y  + y \\
&= 0.
\end{align*}        

This proves (6) for Case 1.

The same argument works for other cases.

Case 2.  $k_{\ell+1} = 2$ and $pk_{\ell+1} = 4n+2$. 
Then, $y_{2\ell+1} = A_{\ell+1}^p y_{2\ell}A_{\ell+1}^{-p}$.\\
 Since $A_{\ell+1}^p = (R_0^{y u_{2\ell}} yx)^p$ 
and  $y_{2\ell} = y_{2\ell-1} = R_0^{w_{2\ell-1}}y$, we have
\begin{equation}
y_{2\ell+1}= R_0^{Q_{p-1} y u_{2\ell}}
(yx)^p R_0^{w_{2\ell-1}}y (yx)^{-p} R_0^{-Q_{p-1} y u_{2\ell}}       
=   R_0^\tau y, \ {\rm where}
\end{equation}
\begin{equation}
\tau = (1-y) Q_{p-1}y u_{2\ell} 
+ (yx)^p w_{2\ell-1}     - (yx)^{n+1} + y.
\end{equation}
Since by (4),
$u_{2\ell} = 
\sum_{j=1}^{\ell} 
m_j (x-1) y^{-1} w_{2j-1}$  and  
$(1-y) Q_{p-1} y (x-1) = 0$,  
by (4.5)(2), we have 
 $\tau = (yx)^p w_{2\ell-1} - (yx)^{n+1} + y 
 = - w_{2\ell-1} - (yx)^{n+1} + y$.

Again, we consider four subcases.\\
Case (i) $\sum_{j=1}^{\ell} k_j \equiv 0$  (mod 4), 
and $\sum_{j=1}^{\ell+1} k_j \equiv 2$ (mod 4).\\
Then by induction   $w_{2\ell-1}$ = 0 and hence 
$w_{2\ell+1} =  -(yx)^{n+1} + y$. \\
Case (ii)   $\sum_{j=1}^{\ell} k_j \equiv 1$  (mod 4), 
and $\sum_{j=1}^{\ell+1} k_j \equiv 3$ (mod 4).\\
Then by induction $w_{2\ell-1} = y$ and hence  
$w_{2\ell+1} = -(yx)^{n+1}$. \\
Case (iii)  $\sum_{j=1}^{\ell} k_j \equiv 2$  (mod 4), 
and $\sum_{j=1}^{\ell+1} k_j \equiv 0$ (mod 4).\\
Since $w_{2\ell-1} = y - (yx)^{n+1}, w_{2\ell+1}$ = 0.\\
Case (iv)  $\sum_{j=1}^{\ell} k_j \equiv 3$  (mod 4), 
and $\sum_{j=1}^{\ell+1} k_j \equiv 1$ (mod 4).\\
Since $w_{2\ell-1} =  - (yx)^{n+1}$, 
$w_{2\ell+1} =  y$.   
This proves (6) for Case 2.   

Case 3.  $ k_{\ell+1}= 3$ and $pk_{\ell+1}= 6n+3$.\\
Then, $y_{2\ell+1} = 
A_{\ell+1}^{3n+2} x_{2\ell}A_{\ell+1}^{-(3n+2)}$.
 Since $A_{\ell+1}^{3n+2} = (R_0^{y u_{2\ell}} yx)^{3n+2}\\
= R_0^{Q_{3n+1}yu_{2\ell}}(yx)^{3n+2}$, 
it follows that  $y_{2\ell+1} = R_0^{\tau}y$,  where 
\begin{equation}
\tau = \{(1-y) Q_{3n+1} y + (yx)^{3n+2}\}
u_{2\ell} - (yx)^{3n+2}y^{-1}w_{2\ell-1} - (yx)^{n+1}.
\end{equation}                                
Since $\{(1-y) Q_{3n+1}y + (yx)^{3n+2}\}(x-1) = 0$  
by (4.5)(3), we have \\
 $\tau  = - (yx)^{3n+2}y^{-1}w_{2\ell-1} 
 - (yx)^{n+1} = (yx)^{n+1}y^{-1}w_{2\ell-1} 
 - (yx)^{n+1}$ .

\noindent
Case (i)   $\sum_{j=1}^{\ell} k_j \equiv 0$  (mod 4), 
and  $\sum_{j=1}^{\ell+1} k_j \equiv 3$ (mod 4).\\
Then by induction $w_{2\ell-1} = 0$ and 
$w_{2\ell+1} =  -(yx)^{n+1}$. \\
Case (ii)   $\sum_{j=1}^{\ell} k_j \equiv 1$  (mod 4), 
and $\sum_{j=1}^{\ell+1} k_j \equiv 0$ (mod 4).\\
Then $w_{2\ell-1} = y$ and hence $w_{2\ell+1}$ = 0. \\
Case (iii)    $\sum_{j=1}^{\ell} k_j \equiv 2$  (mod 4), 
and $\sum_{j=1}^{\ell+1} k_j \equiv 1$ (mod 4).\\
Since $w_{2\ell-1}= y-(yx)^{n+1}$, we have \\
$w_{2\ell+1}= (yx)^{n+1} -(yx)^{n+1}(xy)^n x - (yx)^{n+1} 
= y$.\\
Case (iv)    $\sum_{j=1}^{\ell} k_j \equiv 3$  (mod 4), and 
$\sum_{j=1}^{\ell+1} k_j \equiv 2$ (mod 4).\\
Since $w_{2\ell-1} =  - (yx)^{n+1}$, $w_{2\ell+1} =  
y - (yx)^{n+1}$.\\
This proves (6) for Case 3, 
and the proof of the first part of Theorem A is complete.

\section{
Proof of Theorem A. (II), Proof of Proposition 5.2(2)}

In this section, we prove that if $K(r)$ is a knot, 
then
\begin{equation}
\widetilde \lambda (r) = 
\sum_{j=1}^m (-1)^{\ell(u_j)}\epsilon_j u_j = \pm \mu w,
\ {\rm where}\ w \in  F(x,y)\ {\rm  and}\ \mu \in \ZZ[s_0].
\end{equation}

For simplicity, to each element $u$ in 
$\widetilde A(s_0)$, say  $u = \sum_{j} \epsilon_j u_j$, 
we write  
$\widetilde u = \sum_{j} (-1)^{\ell(u_j)}\epsilon_j u_j$.

First we notice a similar proposition to Proposition \ref{prop:6.3}
holds.  
Since a proof is exactly the same, 
we omit the details.

\begin{prop}\label{prop:8.1}
Let $r = [pk_1,2m_1, pk_2, 2m_2, \cdots, 2m_q, pk_{q+1}]$ and 
$r' = [pk'_1,2m_1, pk'_2, 2m_2, \cdots, 2m_q, pk'_{q+1}]$.
Then  $\widetilde\lambda (r) = \widetilde\lambda (r')$   
if $ k  \equiv  k'$  (mod 4)  for $j = 1,2,\cdots,2q +1$.
\end{prop}

To evaluate 
$\widetilde w_j, \widetilde u_j$ and $\widetilde v_j$, 
we repeat the same argument that 
was used in Section~7.  
But we employ Lemma \ref{lem:4.5} 
instead of Lemma \ref{lem:4.4}. 

\begin{prop}\label{prop:8.2}
Let $r = [pk_1,2m_1, pk_2, 2m_2, 
\cdots, 2m_q, pk_{q+1}]$, where $k_j=1,2$, or $3$ 
for any $j \geq 1$.  
%
%
%
%
%
%
%
%
%
 Then we have
\begin{align}
&(1)\ \widetilde w_0 = 0, 
 \widetilde u_0 = 0\ {\it  and}\ 
 \widetilde v_0 = 0.\nonumber\\
&(2)\  {\it For}\ j \geq 1,   
\widetilde w_{2j-1} = \widetilde w_{2j}, \ {\it and}\   
\widetilde u_{2j-2} = \widetilde u_{2j-1}.\nonumber\\
&(3)\  {\it For\ any}\ 
j \geq 0, 
\widetilde v_j =\widetilde u_j + y^{-1} \widetilde w_j.
\nonumber\\
&(4)\ {\it Suppose}\  q \geq 0.\nonumber\\
%
&\ \ (a)\ {\it If}\ k_{q+1} = 1, \ {\it then}\nonumber\\
&\ \ \ \ \ \widetilde w_{2q+1} = \{-(1+y) Q_n y +(yx)^{n+1}\}  
\widetilde u_{2q} +(yx)^{n+1}y^{-1}\widetilde w_{2q}-y.\nonumber\\
&\ \ (b)\ {\it If}\  k_{q+1} = 2,\ {\it  then}\nonumber\\      
&\ \ \ \ \  \widetilde w_{2q+1} = \{-(1+y) Q_{2n}y\} 
\widetilde u_{2q} -\widetilde w_{2q} -(yx)^{n+1}-y.\nonumber\\
&\ \ (c)\ {\it If}\  k_{q+1} = 3,\ {\it  then}\nonumber\\
&\ \ \ \ \ 
\widetilde w_{2q+1} = \{-(1+y)Q_{3n+1}y - (yx)^{n+1}\} 
\widetilde u_{2q} - (yx)^{n+1}y^{-1} \widetilde w_{2q} -(yx)^{n+1}.\nonumber\\
%
&(5)\ {\it For\  any}\  
j \geq 1, \widetilde u_{2j} = m_j (x+1) y^{-1} \widetilde w_{2j-1} + \widetilde u_{2j-1}.
\end{align}       

\end{prop}

\noindent{\it Proof.}  (1) and (2) follow immediately, noting  
$\widetilde x   = - x$.  Also 
(3) follows from Proposition \ref{prop:7.1}(5), 
since  $\widetilde y  = -y$. 
Next we prove (5).   
Consider $B_j = x_{2j-1} z_{2j-1}^{-1} =  
R_0^{v_{2j-1}-u_{2j-1}} =  R_0^{-y^{-1}w_{2j-1}}$.
Then by (7.3), we see 
$z_{2j} = B_j^{m_j} z_{2j-1}B_j^{-m_j}       
 =  R_0^{-m_j (1-x) y^{-1}w_{2j-1} + u_{2j-1}} x$,
and hence,
$u_{2j} = - m_j (1-x) y^{-1} w_{2j-1} + u_{2j-1}$ and 
$\widetilde u_{2j} =  m_j (1+x) y^{-1} 
\widetilde w_{2j-1} + \widetilde u_{2j-1}$.

Finally, we prove (4) by induction.
For the initial case $q=0$, Proposition \ref{prop:8.2} 
holds.  In fact, if $k_1 = 1$, (7.4) shows that $w_1 = y$ 
and $\widetilde w_1 = -y$.   
If $k_1  = 2$, then from (7.5) we 
see that $w_1 = -(yx)^{n+1} + y$ and 
$\widetilde w_1 =  -(yx)^{n+1}-y$.  If $k_1 = 3$, then 
$w_1 = -(yx)^{n+1} = \widetilde w_1$ by (7.6).

Now suppose  Proposition \ref{prop:8.2}(4) 
holds for $q$ and prove it for $q+1$.\\
If $k_{q+1} =1$, then (7.7) yields, since 
$v_{2q} = u_{2q}-y^{-1} w_{2q}$, 
$w_{2q+1} = (1-y) Q_n  y u_{2q} +(yx)^{n+1}(u_{2q}
-y^{-1}w_{2q}) + y$.\\
By taking a tilde on each element in both sides, 
we obtain (4)(a).
If $k_{q+1}= 2$, 
then since $(yx)^{2n+1} = -1$, 
$w_{2q+1}= (1-y) Q_{2n}y u_{2q}-w_{2q} -(yx)^{n+1} + y$
by (7.10).
By taking a tilde on each element, we have (4)(b).\\
  If $k_{q+1}= 3$, then (7.11) yields\\            
$w_{2q+1} = \{(1-y)Q_{3n+1}y + (yx)^{3n+2}\}u_{2q}- (yx)^{3n+2} 
y^{-1}w_{2q} -(yx)^{n+1}$
and since $(yx)^{2n+1}=1$,
(4)(c) follows by taking a tilde on each element. 
This proves Proposition \ref{prop:8.2}. 
\qed

\begin{thm}\label{thm:8.3}
Let $r = 
[pk_1,2m_1, pk_2, 2m_2, \cdots, 2m_q, pk_{q+1}]$, 

 $r' = [pk_1,2m_1, pk_2, 2m_2, \cdots, 2m_{q-1}, pk_q]$, and
 
$\widehat r=  
[pk_1,2m_1, pk_2, 2m_2, \cdots, 2m_{q-1}, p(k_q + k_{q+1})]$. 

Then $\widetilde\lambda (r)$ is of the form:
\\
(1) If
$\sum_{j=1}^{q+1} k_j \equiv 0\ ({\rm mod}\ 4), \ then\  
\widetilde\lambda (r) = F_0(r) (y- (yx)^{n+1})$.\\  
(2) If  
$\sum_{j=1}^{q+1} k_j \equiv 1\ ({\rm mod}\ 4),\ then\   
\widetilde\lambda (r) = F_1(r) y$.\\
(3) If
$\sum_{j=1}^{q+1} k_j \equiv 2\  ({\rm mod}\  4),\  then\   
\widetilde\lambda (r) = F_2(r) (y+ (yx)^{n+1})$.\\  
(4) If 
$\sum_{j=1}^{q+1} k_j \equiv 3\ ({\rm mod}\ 4),\  then\   
\widetilde\lambda (r) = F_3(r) (yx)^{n+1}$.
\hfill \mbox{\rm{(8.3)}}

\setcounter{equation}{3}

Here $F_j(r),0 \leq  j \leq 3$, 
are complex numbers in $\ZZ[s_0]$, and these 
numbers are determined inductively as follows. \\
(I)
$F_0([0]) = 0,  F_1([p]) =  F_2([2p]) = F_3([3p]) = -1$\\ 
(II) Suppose  $q > 0$.
\begin{align}
&
 (1)\ {\it If}\  k_{q+1}=1,\  {\it then}\ \nonumber\\
&\ \ 
(i)\  F_0(r) = 4m_q b_n  F_3 (r') + F_0(\widehat r)\nonumber\\
&\ \ 
(ii)\   F_1(r) = -8m_q b_n  F_0 (r') + F_1(\widehat r)\nonumber\\
&\ \ 
 (iii)\   F_2(r) = -4m_q b_n  F_1 (r') + F_2(\widehat r)\nonumber\\
&\ \ 
(iv)\  F_3(r) = -8m_q b_n  F_2 (r') + F_3(\widehat r)\nonumber
\end{align}
\begin{align}
&
 (2)\ {\it If}\  k_{q+1} = 2,\ {\it then}\nonumber\\
&\ \ 
 (i)\ F_0(r) = 8m_q b_n  F_2 (r') + F_0(\widehat r)\nonumber\\
&\ \ 
 (ii)\   F_1(r) = 8m_q b_n  F_3 (r') + F_1(\widehat r)\nonumber\\
&\ \ 
(iii)\  F_2(r) = -8m_q b_n  F_0(r') + F_2(\widehat r)\nonumber\\
&\ \ 
(iv)\  F_3(r) = -8m_q b_n  F_1(r') + F_3(\widehat r)\nonumber
\end{align}
\begin{align}
&
(3)\ {\it If}\ k_{q+1} = 3,\ {\it then}\nonumber\\
&\ \ 
(i)\  F_0(r) = 4m_q b_n  F_1(r') + F_0(\widehat r)\nonumber\\
&\ \  
(ii)\ F_1(r) = 8m_q b_n  F_2 (r') + F_1(\widehat r)\nonumber\\
&\ \  
(iii)\ F_2(r) = 4m_q b_n  F_3 (r') + F_2(\widehat r)\nonumber\\
&\ \ 
(iv)\ F_3(r) = -8m_q b_n  F_0(r') + F_3(\widehat r)
\end{align}

Here $b_n$ is the $(1,2)$ entry of the matrix $(XY)^n$, 
see (4.1).
\end{thm}

\begin{rem}\label{rem:8.1}
We use these formulas 
as follows.
For example, 
suppose $k_{q+1}=1$.  
If   
$\sum_{j=1}^{q+1} k_j \equiv 0$ (mod 4), 
then we see by (8.3)(1), 
$\widetilde \lambda (r) = F_0(r) (y - (yx)^{n+1})$.  
In this case, since  $k_{q+1}=1$, 
it follows $\sum_{j=1}^q k_j \equiv 3$ (mod 4)  
and hence by (8.3)(4), we see 
 $\widetilde\lambda (r') = F_3(r') (yx)^{n+1}$.       
Further, 
$\sum_{j=1}^{q -1} k_j  + (k_q + k_{q+1}) 
\equiv 0$ (mod 4) implies that
$\widetilde \lambda (\widehat r) = 
F_0(\widehat r) (y-(yx)^{n+1})$.    
We know inductively $F_3(r')$ and 
$F_0(\widehat r)$, since the lengths of $r'$ and 
$\widehat r$ are shorter than that of $r$, 
and therefore, $F_0(r)$ is determined by (8.4)(II)
(1)(i)
using  $F_3(r')$ and  $F_0(\widehat r)$.  
We list $F_j(r)$ for $q=1$ in the next section.
\end{rem}

\noindent
{\it Proof of Theorem \ref{thm:8.3}.}
We use induction on $q$.  For the initial 
case, $q = 0$, since  
$r = [pk_1]$ and 
$r'= \widehat r = [0]$, (8.3) and (8.4)(I) 
follow from (8.2)(4).  
Note $\widetilde u_0 = \widetilde w_0=0$. 

Next consider the case $q = 1$, i.e. 
$r = [pk_1, 2m_1, pk_2]$.  Then \\
$r'=[pk_1]$ and $\widehat r = [p(k_1 + k_2)]$. 
  
Case (1) $k_2=1$. \\
By (8.2)(4)(a), we  have 
 $\widetilde w_3 = [-(1+y) Q_n y +(yx)^{n+1}]
 \widetilde u_2 +(yx)^{n+1}y^{-1}\widetilde w_2 - y$.\\
Since $\widetilde u_2 = m_1(x+1) y^{-1}
\widetilde w_1$ and  
$\widetilde w_2=\widetilde w_1$,    
we see\\
$\widetilde w_3 = [-(1+y) Q_n y +(yx)^{n+1}] m_1(1+x) y^{-1}   \widetilde w_1 +(yx)^{n+1}y^{-1} \widetilde w_1- y$.\\
Further by (4.6)(1), we have 
\begin{equation}
\widetilde w_3 =  -4b_n m_1( y +(yx)^{n+1}) y^{-1}  
\widetilde w_1 +(yx)^{n+1}y^{-1}\widetilde w_1 - y.
\end{equation}

Now we apply (8.3).  \\
If $k_1=1$, then $\widetilde w_1 = -y$, and hence 
$\widetilde w_3 = (4b_n m_1 - 1) (y +(yx)^{n+1})$.
This proves (8.3) for this case.\\
If $k_1=2$, then  
$\widetilde w_1=  - (y +(yx)^{n+1})$, and hence\\
$\widetilde w_3 = 4b_n m_1( y +(yx)^{n+1}) 
y^{-1}(y +(yx)^{n+1}) 
- (yx)^{n+1}y^{-1}(y + (yx)^{n+1}) - y\\
\hspace*{0.45cm}
= (8b_n m_1 - 1) (yx)^{n+1}$.\\
If  $k_1=3$, then  $\widetilde w_1 =  - (yx)^{n+1}$, 
and hence\\
$\widetilde w_3 =  
4b_n m_1(y +(yx)^{n+1})y^{-1}(yx)^{n+1}
- (yx)^{n+1}y^{-1}(yx)^{n+1} - y\\
\hspace*{0.45cm}
 = 4b_n m_1((yx)^{n+1} - y)$.\\
This proves (8.3) for Case (1), $k_2=1$.

Since similar arguments work for other cases, 
we skip details.
\medskip

Case (2)  $k_2=2$.\\
By (8.2)(4)(b), we  have 
$\widetilde w_3 = \{-(1+y) Q_{2n}y\} 
\widetilde u_2 - \widetilde w_2 - (yx)^{n+1} - y$.\\
By (4.6)(3), it becomes to 
\begin{equation}
\widetilde w_3 = -8b_n m_1(yx)^{n+1}y^{-1} 
\widetilde w_1 -\widetilde w_1 - (yx)^{n+1} - y.
\end{equation}
As before, compute 
$\widetilde w_3$  to each case 
$k_1=1,2$ or $3$ to prove (8.3).
 
Case (3) $k_2=3$.\\
By (8.2)(4)(c), we  see \\
 $\widetilde w_3 = \{-(1+y) Q_{3n+1}y  - (yx)^{n+1}\} 
 \widetilde u_2 - (yx)^{n+1}y^{-1} 
 \widetilde w_2 - (yx)^{n+1}$. \\
By (4.6)(5), it becomes to 
\begin{equation}
\widetilde w_3 =  4b_n m_1 (y - (yx)^{n+1}) y^{-1}  
\widetilde w_1 - (yx)^{n+1}y^{-1}\widetilde w_1 - (yx)^{n+1}.
\end{equation}      
Computation of  $\widetilde w_3$  to each case 
$k_1=1,2$ or $3$ completes the 
proof for $q=1$.

\medskip
Next we assume that Theorem \ref{thm:8.3} 
holds for any $r$ with length less 
than $2q+1$.
First consider the case where $k_{q+1}=1$.   
We divide our proof into three subcases.

Case (1.1)   $(k_q, k_{q+1}) = (1,1)$.\
From (8.2)(4)(a), we have\\
$\widetilde w_{2q+1} = \{-(1+y) Q_n y +(yx)^{n+1}\}  
\widetilde u_{2q} +(yx)^{n+1}y^{-1}\widetilde w_{2q} - y$.\\
Since $\widetilde u_{2q} = m_q (1+x) y^{-1}
\widetilde w_{2q-1} + \widetilde u_{2q-1}$, it follows that
\begin{align*}
\widetilde w_{2q+1}
&= 
\{-(1+y) Q_n y +(yx)^{n+1}\} (m_q (1+x) y^{-1}
\widetilde w_{2q-1} + \widetilde u_{2q-1})\\
&\ \  +(yx)^{n+1}y^{-1}\widetilde w_{2q}~-~y\\
&=  \{-(1+y) Q_n y +(yx)^{n+1}\} m_q (1+x) y^{-1}
\widetilde w_{2q-1} \\
&\ \ + \{-(1+y) Q_n y +(yx)^{n+1}\} 
\widetilde u_{2q-2}\\
&\ \ +(yx)^{n+1}y^{-1}\widetilde w_{2q-1} - y,
\end{align*}
since $\widetilde w_{2q}= \widetilde w_{2q-1}$ and  $\widetilde u_{2q-1}= \widetilde u_{2q-2}$.

Let $A =  \{-(1+y) Q_n y +(yx)^{n+1}\} (1+x)$  and \\
    $B = \{-(1+y) Q_n y +(yx)^{n+1}\} 
    \widetilde u_{2q-2} +(yx)^{n+1}y^{-1}\widetilde w_{2q-1} - y$.\\
Then  $\widetilde w_{2q+1} = A m_q y^{-1} 
\widetilde w_{2q-1} + B$.

First we claim that $B = \widetilde \lambda (\widehat r)$.
To prove this claim we should note that, since 
$k_q=1$, by induction assumption,  \\
$\widetilde w_{2q-1} = \{-(1+y) Q_n y +(yx)^{n+1}\}  
\widetilde u_{2q-2} +(yx)^{n+1}y^{-1}\widetilde w_{2q-2} -~y$.
Therefore, 
\begin{align*}
B &= 
\{-(1+y) Q_n y +(yx)^{n+1}\} 
\widetilde u_{2q-2} \\
&\ \ +(yx)^{n+1}y^{-1}
\Bigl[\{-(1+y) Q_n y +(yx)^{n+1}\}
 \widetilde u_{2q-2} 
  +(yx)^{n+1}y^{-1}\widetilde w_{2q-2}-y \Bigr] - y\\
&=\Bigl[-(1+y) Q_n y +(yx)^{n+1} +(yx)^{n+1}y^{-1} \bigl\{-(1+y) 
Q_n y +(yx)^{n+1} \bigr\}\Bigr] \widetilde u_{2q-2} \\
&\ \ +(yx)^{n+1}
y^{-1}\{(yx)^{n+1} y^{-1}\widetilde w_{2q-2}-y\} - y.
\end{align*}
 
Since $k_q + k_{q+1}=2$, 
it suffices to show, using (8,2)(4)(b),
\begin{align}
& (i)\   
-(1+y) Q_n y +(yx)^{n+1} - (yx)^{n+1}y^{-1}(1+y) 
Q_n y +(yx)^{n+1}y^{-1}(yx)^{n+1}\nonumber\\
&\ \ \  \ =-(1+y) Q_{2n}y,\  {\rm and}\nonumber\\ 
& (ii)\ 
(yx)^{n+1}y^{-1}(yx)^{n+1}y^{-1}=-1.
\end{align}          

\noindent
{\it Proof of (8.8).} (ii) follows immediately, 
and then, (i) becomes to\\
$-(1+y) Q_n y +(yx)^{n+1} -(yx)^{n+1}y^{-1}(1+y) 
Q_n y - y =  -(1+y)Q_{2n}y$.\\
Since
$(yx)^{n+1} = (yx)(yx)^n = y(yx)^n y$, 
the above equation is  
equivalent to  \\
(i)'
$-(1+y)Q_n + y(yx)^{n} - y(yx)^n (1+y)Q_n -1  
= -(1+y) Q_{2n}$.               \\
Since $Q_{2n} = Q_n + (yx)^n (Q_n - 1)$, we see\\
LHS (of (i)') \\
$=-(1+y) Q_{2n} + (1+y) (yx)^n (Q_n -1) + y(yx)^n     
- y(yx)^n (1+y) Q_n -1\\
=-(1+y) Q_{2n} + (1+y) (yx)^n Q_n - (1+y) (yx)^n + y(yx)^n     
- y(yx)^n (1+y) Q_n -1\\
=-(1+y)Q_{2n} + \Bigl\{(1+y) (yx)^n - y(yx)^n (1+y)
\Bigr\}Q_n - (yx)^n -1\\
=-(1+y) Q_{2n} + (yx)^n (1 - yx) Q_n - (yx)^n - 1\\
=-(1+y) Q_{2n} + (yx)^n (1 - (yx)^{n+1}) - (yx)^n - 1\\
=-(1+y) Q_{2n} - (yx)^{2n+1} - 1\\
=- (1+y) Q_{2n}$. \\
This proves (8.8) and $B=\widetilde\lambda(\widehat{r})$.
Therefore, we have 
\begin{equation}
\widetilde w_{2q+1} = A m_q y^{-1}\widetilde w_{2q-1}
 + \widetilde \lambda (\widehat r).
 \end{equation}    

We note (4.6)(1) shows us that 
$A = -4b_n (y + (yx)^{n+1})$.

To prove (8.4)(II)(1), 
we consider the following four cases separately.

Case (i)   Suppose  $\sum_{j=1}^{q+1}k_j  \equiv 0$  (mod 4).
Since $k_{q+1}=1, \sum_{j=1}^q k_j  \equiv  3$ (mod 4) 
and hence, by induction 
assumption, $\widetilde w_{2q-1} = F_3(r') (yx)^{n+1}$.  
Therefore
\begin{align*}
A m_q y^{-1} \widetilde w_{2q-1} 
&= m_q F_3 (r') \Bigl\{- 4b_n (y +(yx)^{n+1})
\Bigr\}y^{-1}(yx)^{n+1}\\       
&= -4 m_q F_3(r') b_n ((yx)^{n+1} -y).
\end{align*}
Also, by induction,  
$\widetilde\lambda (\widehat r) = F_0(\widehat r) 
(y - (yx)^{n+1})$, and hence \\
$\widetilde\lambda (r) =  
(4 m_q F_3(r') b_n + F_0(\widehat r)
) (y-(yx)^{n+1})$.
This proves (8.4)(II)(1)(i).

Case (ii)  $\sum_{j=1}^{q+1}k_j  \equiv 1$  (mod 4).
Then $\sum_{j=1}^q k_j  \equiv  0$ (mod 4) and 
 hence, by induction assumption,
we obtain that $\widetilde w_{2q-1}= F_0(r') 
(y - (yx)^{n+1})$, and therefore, 
\begin{align*}
A m_q y^{-1}\widetilde w_{2q-1}
& = -4 m_q F_0(r') b_n
(y +(yx)^{n+1}) y^{-1}(y-(yx)^{n+1})\\
&=- 4 m_q F_0(r') b_n y(1 +(xy)^n x) (1 - (xy)^n x)\\
&=-4 m_q F_0(r') b_n y (1 - (xy)^n xy (xy)^n)\\
&=- 8 m_q F_0(r') b_n y.
\end{align*}
Also, by induction,  
$\widetilde\lambda (\widehat r) =
 F_1(\widehat r) y$, and hence 
$\widetilde\lambda (r) =  
(-8 m_q F_0(r') b_n + F_1   (\widehat r)) y$.

Case (iii).   $\sum_{j=1}^{q+1}k_j  \equiv 2$  (mod 4).
Then $\sum_{j=1}^q k_j  \equiv  1$ (mod 4) and 
 hence, by induction 
assumption,  $\widetilde\lambda (r') = F_1(r')y$, and 
\begin{align*}
A m_q y^{-1} \widetilde w_{2q-1} 
&= -4 m_q F_1(r') b_n   
(y +(yx)^{n+1}) y^{-1}y\\
&=  - 4 m_q F_1(r') b_n (y +(yx)^{n+1}).
\end{align*}
                     
On the other hand, $\widetilde\lambda (\widehat r) = 
F_2(\widehat r) (y + (yx)^{n+1})$, and hence \\
$\widetilde\lambda (r) = (-4 m_q F_1(r') 
b_n + F_2(\widehat r)) (y +(yx)^{n+1})$.

Case (iv)     $\sum_{j=1}^{q+1}k_j  \equiv 3$  (mod 4).
Then $\sum_{j=1}^q k_j  \equiv 2$ (mod 4) and 
hence, $\widetilde\lambda (r') = F_2(r') (y + (yx)^{n+1})$, and
\begin{align*}
A m_q y^{-1} \widetilde w_{2q-1} 
&= 
-4 m_q F_2(r')b_n 
(y +(yx)^{n+1})y^{-1}(y + (yx)^{n+1})\\
&=  - 4 m_q F_2(r')b_n y(1+ (xy)^n x) (1+ (xy)^n x)\\
&=  - 4 m_q F_2(r')b_n y (1+ 2(xy)^n x+ (xy)^{2n+1})\\
&=   -8 m_q F_2(r') b_n (yx)^{n+1}.
\end{align*}                          

Also, by induction,  
$\widetilde\lambda (\widehat r) = F_3(\widehat r) (yx)^{n+1}$, 
and hence \\
   $\widetilde \lambda (r) = (-8 m_q F_2(r') 
   b_n + F_3(\widehat r)) (yx)^{n+1}$. \\
Therefore,  Theorem \ref{thm:8.3} is proved for this case.

For other cases, we use essentially the same argument, 
although  calculations for some cases are a bit complicated. 
We just state the final forms and details will be omitted.

Case (2.1)    $(k_q, k_{q+1}) = (2,1)$

First we write  $\widetilde w_{2q+1}= A m_q y^{-1} 
\widetilde w_{2q-1} + B$, where \\
$A = \{-(1+y) Q_n y +(yx)^{n+1}\} (1+x)$  and\\
$B =  \{-(1+y)Q_n y +(yx)^{n+1} 
-(yx)^{n+1}y^{-1}(1+y)Q_{2n}y\} 
\widetilde u_{2q-2}\\            
\hspace*{5mm} + (yx)^{n+1}y^{-1}( - \widetilde w_{2q-2} - (yx)^{n+1} - y)-y$.\\
Then, we can show that 
\begin{align}
&(1)\ -(1+y) Q_n y +(yx)^{n+1}
-(yx)^{n+1}y^{-1}(1+y) Q_{2n}y\nonumber\\ 
&\ \ \ \ = - (1+y)Q_{3n+1}y -(yx)^{n+1},\ {\rm  and} \nonumber\\
&(2)\  (yx)^{n+1}y^{-1}
(-\widetilde w_{2q-2} - (yx)^{n+1} - y) - y\nonumber\\
 &\ \ \ \ =  - (yx)^{n+1}y^{-1} \widetilde w_{2q-2} - (yx)^{n+1}. 
 \end{align}
Therefore, we see that 
$B = \widetilde\lambda (\widehat r)$, and further,
$A = -4b_n (y + (yx)^{n+1})$, 
and thus,
$\widetilde w_{2q+1} = -4m_q b_n (y + (yx)^{n+1}) y^{-1}  
 \widetilde w_{2q-1} + \widetilde \lambda (\widehat r)$.\\

Case (i)   Suppose   $\sum_{j=1}^{q+1}k_j  \equiv 0$  (mod 4).\\
Then $\sum_{j=1}^q k_j  \equiv 3$ (mod 4) and hence, 
 by induction assumption,  
$\widetilde w_{2q-1} = F_3(r') (yx)^{n+1}$, and 
\begin{align*}
A m_q y^{-1}\widetilde w_{2q-1}
&= - 4m_q F_3(r') b_n   
(y + (yx)^{n+1}) y^{-1}(yx)^{n+1} \\                   
&= -4 m_q F_3(r') b_n ((yx)^{n+1} -y).
\end{align*}
Therefore,
$\widetilde \lambda (r) = (4m_q F_3(r')b_n + 
F_0(\widehat r)) (y-(yx)^{n+1})$.

Case (ii).    $\sum_{j=1}^{q+1}k_j  \equiv 1$  (mod 4).\\
Then $\sum_{j=1}^q k_j  \equiv  0$ (mod 4) and hence, 
$\widetilde w_{2q-1} = F_0(r') (y - (yx)^{n+1})$, and  
\begin{align*}
A m_q y^{-1}\widetilde w_{2q-1} 
&= -4 m_q F_0(r') 
b_n   (y +(yx)^{n+1}) y^{-1}(y-(yx)^{n+1})\\                     
&= - 8 m_q F_0(r') b_n y. 
\end{align*}
Also by induction, 
$\widetilde \lambda (\widehat r) = F_1(\widehat r) y$, 
and hence 
    $\widetilde \lambda (r) = (-8 m_q F_0(r') 
    b_n + F_1   (\widehat r)) y$.
    
Case (iii)     $\sum_{j=1}^{q+1}k_j  \equiv 2$  (mod 4).\\
Then $\sum_{j=1}^q k_j  \equiv 1$ (mod 4) and, 
by induction assumption,  
$\widetilde w_{2q-1} = F_1(r')y$.    
Therefore, 
$A m_q y^{-1}\widetilde w_{2q-1} = -4 m_q F_1(r') b_n   
(y +(yx)^{n+1}) y^{-1}y$.                        

On the other hand, 
$\widetilde\lambda (\widehat r) 
= F_2(\widehat r) (y + (yx)^{n+1})$, and hence \\
   $\widetilde\lambda (r) = (-4 m_q F_1(r') b_n 
   + F_2(\widehat r)) (y +(yx)^{n+1})$.

Case (iv)     $\sum_{j=1}^{q+1}k_j  \equiv 3$  (mod 4).\\
Then $\sum_{j=1}^q k_j  \equiv  2$ (mod 4) 
and, by induction 
assumption, $\widetilde w_{2q-1} = F_2(r') 
(y + (yx)^{n+1})$.  Therefore, 
\begin{align*}
A m_q y^{-1}\widetilde w_{2q-1} 
&= -4 m_q F_2(r') b_n  (y +(yx)^{n+1})y^{-1}(y + (yx)^{n+1})  \\                   
&=  -8 m_q F_2(r') b_n (yx)^{n+1}.
\end{align*}                      
Since
$\widetilde \lambda (\widehat r) 
= F_3(\widehat r) (yx)^{n+1}$, we have
$\widetilde \lambda (r) =  
(-8 m_q F_2(r') b_n + F_3 (\widehat r)) (yx)^{n+1}$.       
Therefore, for this case, Theorem \ref{thm:8.3} is proved.

Case (3.1)   $(k_q, k_{q+1})=(3,1)$

As above, we write  
$\widetilde w_{2q+1} =  A m_q y^{-1}  
\widetilde w_{2q-1} + B$,
where\\
 $A = \{-(1+y) Q_n y +(yx)^{n+1}\} (1+x) 
= -4b_n (y + (yx)^{n+1})$, and \\
$B = \{-(1+y) Q_n y +(yx)^{n+1}+(yx)^{n+1}y^{-1}
(- (1+y) Q_{3n+1}y - (yx)^{n+1})\} \widetilde u_{2q-2} \\  
\hspace*{5mm}+ (yx)^{n+1} y^{-1}( - (yx)^{n+1}y^{-1} 
\widetilde w_{2q-2} - (yx)^{n+1})- y$.
We can show that $B = \widetilde\lambda (\widehat r)$.  
 Therefore,
$\widetilde w_{2q+1} = -4m_q b_n (y + (yx)^{n+1}) y^{-1}  
\widetilde  w_{2q-1} + \widetilde \lambda (\widehat r)$. 

Case (i)   Suppose   $\sum_{j=1}^{q+1}k_j  \equiv 0$  (mod 4).\\
Then $\sum_{j=1}^q k_j  \equiv 3$ (mod 4) and hence, 
\begin{align*}
A m_q y^{-1}\widetilde w_{2q-1} 
&= -4 m_q F_3(r') b_n   (y +(yx)^{n+1}) y^{-1}(yx)^{n+1}\\
&= -4 m_q F_3  (r') b_n ((yx)^{n+1} -y).
\end{align*}
And thus,
$\widetilde \lambda (r) = 
(4m_q F_3(r') b_n + F_0(\widehat r)) (y - (yx)^{n+1})$.

Case (ii)     $\sum_{j=1}^{q+1}k_j  \equiv 1$  (mod 4).\\
Then $\sum_{j=1}^q k_j  \equiv  0$ (mod 4) and, since 
$\widetilde w_{2q-1}= F_0(r') (y - (yx)^{n+1})$, 
\begin{align*}
A m_q y^{-1}\widetilde w_{2q-1} 
&= -4m_q F_0(r') 
b_n  (y +(yx)^{n+1}) y^{-1}(y-(yx)^{n+1})  \\                   
&=  - 8 m_q F_0(r') b_n y.
\end{align*}                     
Also by induction, 
$\widetilde\lambda (\widehat r) = F_1(\widehat r) y$, and hence 
$\widetilde\lambda (r) = (-8 m_q F_0(r') b_n + F_1(\widehat r)) y$.

Case (iii)        $\sum_{j=1}^{q+1}k_j  \equiv 2$  (mod 4).\\
Then $\sum_{j=1}^q k_j  \equiv  1$ (mod 4) and hence,  
$\widetilde w_{2q-1} = F_1(r') y$.\\
Therefore, 
$A m_q y^{-1}\widetilde w_{2q-1} = -4 m_q F_1(r') b_n  
(y +(yx)^{n+1}) y^{-1}y$.                        
On the other hand,   $\widetilde \lambda (\widehat r) = 
F_2(\widehat r) (y + (yx)^{n+1})$, and
$\widetilde \lambda (r) =  
   (-4 m_q F_1(r') b_n + F_2(\widehat r)) (y +(yx)^{n+1})$.

Case (iv)      
$\sum_{j=1}^{q+1}k_j  \equiv 3$  (mod 4).\\
Then $\sum_{j=1}^q k_j  \equiv  2$ (mod 4) and thus,
  $\widetilde w_{2q-1} = F_2(r') (y + (yx)^{n+1})$ and 
\begin{align*}
A m_q y^{-1}\widetilde w_{2q-1} 
&= -4 m_q F_2(r') b_n  
(y +(yx)^{n+1})y^{-1}(y + (yx)^{n+1})\\
&= -8 m_q F_2(r') b_n (yx)^{n+1}.
\end{align*}                          
By induction, since
$\widetilde \lambda (\widehat r) = 
F_3(\widehat r) (yx)^{n+1}$, we have
$\widetilde\lambda (r) = (-8 m_q F_2(r') b_n + 
F_3(\widehat r)) (yx)^{n+1}$.       
For this case, Theorem \ref{thm:8.3} is now proved.\\
From the above proof, we notice that 
$\widetilde \lambda (r)$ depends only on $k_{q+1}$         
and $\sum_{j=1}^{q+1}k_j$  (mod 4).  
Therefore, in the rest of our proof, it 
suffices to consider only the case where 
$(k_q, k_{q+1}) = (1,2)$ and $(1,3)$.

Case (1.2)   $(k_q, k_{q+1})=(1,2)$

We write 
$\widetilde w_{2q+1}= A m_q y^{-1}\widetilde w_{2q-1} + B$,
where $A = - (1+y) Q_{2n} y (1+x)$, and 
$B = \{- (1+y) Q_{2n}y + (1+y)Q_n y - (yx)^{n+1}\}
\widetilde u_{2q-2}          
- (yx)^{n+1} y^{-1}\widetilde w_{2q-2} - (yx)^{n+1}$.      
It is shown that
$B = \{- (1+y) Q_{3n+1} y - (yx)^{n+1}\}
\widetilde u_{2q-2} - (yx)^{n+1}y^{-1}\widetilde w_{2q-2}   
 - (yx)^{n+1}$,      
which is  $\widetilde \lambda (\widehat r)$.  
Further, by (4.6)(3), we see
\begin{equation*}
A m_q y^{-1}\widetilde w_{2q-1} 
= - 8b_n m_q (yx)^{n+1}y^{-1}\widetilde w_ {2q-1}.
\end{equation*}    

 Case (i)   
 $\sum_{j=1}^{q+1}k_j  \equiv 0$  (mod 4).\\
Since $k_{q+1}=2, \sum_{j=1}^q k_j  \equiv  2$ (mod 4) and hence, 
by induction assumption,  
$\widetilde w_{2q-1} = F_2(r') (y +(yx)^{n+1} )$.
Therefore,
\begin{align*}
\widetilde \lambda (r)  &= - 8 m_q F_2(r') b_n 
(yx)^{n+1}  y^{-1}(y +(yx)^{n+1})  
 + F_0(\widehat r) (y - (yx)^{n+1})  \\                  
& = - 8 m_q F_2(r') b_n ((yx)^{n+1} -  y) 
 + F_0(\widehat r) (y - (yx)^{n+1})  \\                  
&=  (8 m_q F_2(r') b_n + F_0(\widehat r)) (y - (yx)^{n+1}).
\end{align*}

Case (ii)
$\sum_{j=1}^{q+1}k_j  \equiv 1$  (mod 4).\\
Then  $\sum_{j=1}^q k_j  \equiv  3$ (mod 4) and 
 $\widetilde w_{2q-1}= F_3(r') (yx)^{n+1}$. 
 Therefore 
 \begin{align*}
 \widetilde\lambda (r)  
 &=  - 8 m_n F_3(r') 
 b_n (yx)^{n+1}y^{-1}(yx)^{n+1} + F_1(\widehat r) y  \\                   
&=  (8 m_n F_3(r') b_n + F_1(\widehat r)) y.
\end{align*}

Case (iii)   
$\sum_{j=1}^{q+1}k_j  \equiv 2$  (mod 4).\\
Then  $\sum_{j=1}^q k_j  \equiv  0$ (mod 4) and hence,               
$\widetilde w_{2q-1} = F_0(r') (y - (yx)^{n+1})$ and
\begin{align*}
\widetilde\lambda (r) 
&= -8 m_q F_0(r') 
b_n (yx)^{n+1} y^{-1}(y - (yx)^{n+1})  
+ F_2(\widehat r) (y + (yx)^{n+1})  \\                      
&= ( -8 m_q F_0(r') b_n + F_2(\widehat r)) (y + (yx)^{n+1}).
\end{align*}                        

Case (iv)     
 $\sum_{j=1}^{q+1}k_j  \equiv 3$  (mod 4).\\
Then  $\sum_{j=1}^q k_j  \equiv  1$ (mod 4) and hence, 
$\widetilde w_{2q-1} = F_1(r') y$ and 
\begin{align*}
\widetilde\lambda (r) 
&=  -8 m_q F_1(r') 
b_n (yx)^{n+1}y^{-1}y + F_3(\widehat r) (yx)^{n+1}  \\
&= (-8 m_q F_1(r') b_n + F_3(\widehat r))(yx)^{n+1}.
\end{align*}  
Thus for this case, Theorem \ref{thm:8.3} is proved.\\

Case (1.3)   $(k_q, k_{q+1}) = (1, 3)$\\
 Let   $\widetilde w_{2q+1} 
 = A m_q y^{-1}\widetilde w_{2q-1} + B$,
where\\
 $A = \{- (1+y)Q_{3n+1} - (yx)^{n+1}\}(x+1)$, and\\ 
$B = \{- (1+y) Q_{3n+1}y - (yx)^{n+1} - 
(yx)^{n+1}y^{-1}{- (1+y) Q_n y 
+ (yx)^{n+1}}\} 
\widetilde u_{2q-2} - (yx)^{n+1}y^{-1}(yx)^{n+1}y^{-1} 
\widetilde w_{2q-2}         
 = \widetilde \lambda (\widehat r)$. \\
Further,
$A m_q y^{-1} \widetilde w_{2q-1} 
= 4b_n m_q (y - (yx)^{n+1}) y^{-1} \widetilde w_{2q-1}$. 
Therefore,
$\widetilde w_{2q+1} =  
4b_n m_q (y - (yx)^{n+1}) y^{-1} 
\widetilde  w_{2q-1} +  \widetilde \lambda (\widehat r)$.
   
Case (i)   
$\sum_{j=1}^{q+1}k_j  \equiv 0$  (mod 4).\\
Then  $\sum_{j=1}^q k_j  \equiv 1$ (mod 4) and hence, 
by induction 
assumption, $\widetilde w_{2q-1} = F_1(r') y$.   
Since $\widetilde \lambda (\widehat r) 
= F_0(\widehat r) (y - (yx)^{n+1})$, 
we have 
\begin{align*}
\widetilde \lambda (r)  
&= 4m_q F_1(r') b_n (y-(yx)^{n+1}) y^{-1}y  +
F_0(\widehat r) (y - (yx)^{n+1})  \\                  
&=  (4 m_q F_1(r') b_n + 
F_0(\widehat r)) (y -(yx)^{n+1}).
\end{align*}

Case (ii)  
 $\sum_{j=1}^{q+1}k_j  \equiv 1$  (mod 4).\\
Since $\sum_{j=1}^q k_j  \equiv  2$ (mod 4), we see 
$\widetilde w_{2q-1} = F_2(r') (y + (yx)^{n+1})$  and  
$\widetilde \lambda ( \widehat r) = F_1(\widehat r) y$, 
and thus
\begin{align*}
\widetilde\lambda (r)  
&=  4 m_q F_2(r') b_n (y - (yx)^{n+1}) y^{-1}(y + (yx)^{n+1}) + 
F_1(\widehat r) y \\                  
&= ( 8 m_q F_2(r') b_n + F_1(\widehat r)) y.
\end{align*} 

Case (iii)
$\sum_{j=1}^{q+1}k_j  \equiv 2$  (mod 4).\\
Then  $\sum_{j=1}^q k_j  \equiv 3$ (mod 4) and hence 
$\widetilde w_{2q-1} = F_3(r') (yx)^{n+1}$ and,   
 $\widetilde \lambda (\widehat r) 
 =  F_2(\widehat r)(y + (yx)^{n+1})$. 
Therefore, 
\begin{align*}
\widetilde \lambda (r) 
&= 4m_n F_3(r') b_n (y - (yx)^{n+1}) y^{-1}(yx)^{n+1}       
 + F_2(\widehat r) (y + (yx)^{n+1})  \\                      
&= ( 4 m_q F_3(r') b_n + F_2(\widehat r)) (y + (yx)^{n+1}).
\end{align*}                        

Case (iv)      $\sum_{j=1}^{q+1}k_j  \equiv 3$  (mod 4).\\
Then  $\sum_{j=1}^q k_j  \equiv  0$ (mod 4) and hence
$\widetilde w_{2q-1} = F_0(r') (y - (yx)^{n+1})$ and   
$\widetilde \lambda (\widehat r) = F_3(\widehat r) (yx)^{n+1}$.
Thus, 
\begin{align*}
\widetilde \lambda (r) 
&=  4 m_q F_0(r') b_n (y - (yx)^{n+1})y^{-1}(y - (yx)^{n+1}) 
+ F_3(\widehat r) (yx)^{n+1}  \\
& = (-8 m_q F_0(r') b_n + F_3(\widehat r)) (yx)^{n+1}.
\end{align*}  
A proof of Theorem \ref{thm:8.3}, and hence,
a proof of Theorem A is now complete.
\qed

\section{Evaluation of $\mu$.}

For $r=[pk_1,2m_1,pk_2,2m_2,\dots,2m_q,pk_{q+1}]$,
we proved that
$\lambda_{\rho,K(r)}(-1)=\mu^2$ for some $\mu\in \ZZ[s_{0}]$.
For convenience, we denote $\mu =\mu(r)$. In this section,
we give an algorithm by which one can compute $\mu(r)$.
We should note that $\mu(r)=F_j(r)$, where $j\equiv \sum_{i=1}^{q+1} k_i$ (mod 4).
As we used in the previous section, let\\
 $r' = [pk_1,2m_1, pk_2, 2m_2,
  \cdots, 2m_{q-1}, pk_q]$, and \\
$\widehat r = [pk_1,2m_1, pk_2, 2m_2, 
\cdots, 2m_{q-1}, p(k_q+k_{q+1})]$.\\
In the proof of Theorem \ref{thm:8.3},
we have shown the following proposition.

\begin{prop}\label{prop:9.1}
The following equalities hold:
\begin{align}
&(1)\ \mu(r)=\nu\mu(r')+\mu(\widehat{r}).\nonumber\\
&(2)\ \mu[0]=0,\mu[p]=\mu[2p]=\mu[3p]=-1.
\end{align}
Here $\nu=m_qb_n\sigma(k_{q+1},M)$ and $b_n$ is the 
$(1,2)$-entry of $\rho(xy)^n$,
$M\equiv \sum_{j=1}^{q+1}k_j$ (mod 4),
$0\le M\le3$, and $\sigma(k_{q+1},M)$ is  given by the following table:
$\begin{array}{llll}
\sigma (1,0) = 4,
&\ \ \ \sigma (1,1) = -8,
&\ \ \ \sigma (1,2) = - 4,
&\ \ \ \sigma (1,3) = - 8,\\
\sigma (2,0) = 8,
&\ \ \ \sigma (2,1) = \ \ 8,
&\ \ \ \sigma (2,2) = - 8,
&\ \ \ \sigma (2,3) = - 8,\\
\sigma (3,0) = 4,
&\ \ \ \sigma (3,1) = \ \ 8, 
&\ \ \ \sigma (3,2) =  \ \ 4,
&\ \ \ \sigma (3,3) = - 8.
\end{array}$

\end{prop}

%
%
%
%


\begin{ex}   
Let $p=3$ and $n=1$, and hence $b_n=1$.

(1)	
Let $r=[3,-4,3,2,3]$. 
To evaluate $\mu(r)$, it is 
convenient to use the tree diagram below:

\begin{figure}[h]
\begin{center}
\includegraphics[height=23mm]{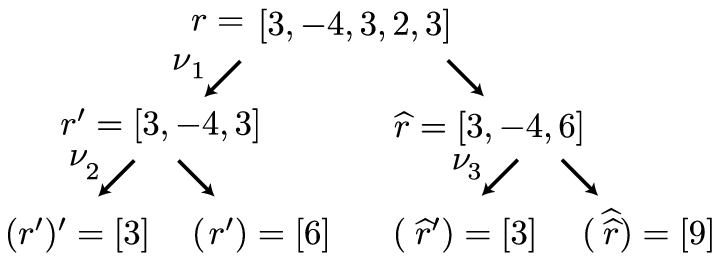}
\end{center}
\end{figure}

%
Since, $m_1 = -2, m_2 = 1$,  
$k_1 = k_2 = k_3 = 1$, 
the weights are
$\nu_1 =\sigma (1,3) = -8,   
\nu_2 = \sigma (1,2) (-2) = (-4)(-2) = 8$
and 
$\nu_3 = \sigma (2,3) (-2) = (-8)(-2) = 16$, 
and hence,
$\mu  = \nu_1 \nu_2 \mu[3] + \nu_1\mu[6] + 
\nu_3\mu[3] + \mu[9] = 55$.
 
(2) 
%
Let $r=[6,2,6,-2,9]$.
Since $m_1 = 1, m_2 = -1,  k_1 = k_2 = 2,  k_3 = 3$, 
the weights are
$\nu_1 = \sigma(3,3)(-1) = (-8)(-1) = 8, \nu_2 = \sigma(2,0)=8$  
and $\nu_3=\sigma(1,3) = -8$, and hence,
$\mu  = \nu_1 \nu_2 \mu[6] + 0 + \nu_3\mu[6] + \mu[9] 
= - 64  + 8 - 1 = - 57$.

\begin{figure}[h]
\begin{center}
\includegraphics[height=23mm]{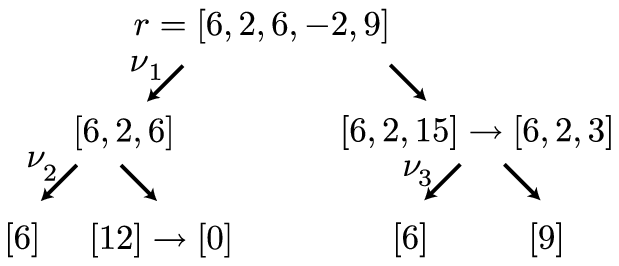}
\end{center}
\end{figure}

\end{ex}

Using these recursion formulas, 
we can prove, for example, the following:
%
%
%
\begin{align}
&(1)\   \mu[p,2m,p] = 4m b_n - 1.\nonumber\\
&(2)\   \mu[p,2m,2p] = 8m b_n - 1.\nonumber\\ 
&(3)\   \mu[p,2m,3p] = - 4m b_n.\nonumber\\
&(4)\   \mu[2p,2m,2p] = -8m b_n.\nonumber\\
&(5)\   \mu[2p,2m,3p] = -8m b_n - 1.\nonumber\\
&(6)\   \mu[3p,2m,3p] = - 4m b_n - 1.
\end{align}

\begin{align}
&(1)\    \mu[p, 2m_1, p, 2m_2,p]
=-32 m_1m_2 b_n^2 + (8m_1 + 8 m_2) b_n - 1\nonumber\\
&(2)\  \mu[p, 2m_1, 2p, 2m_2, 2p] 
= 64 m_1 m_2 b_n^2 - 8m_2 b_n   - 1.
\end{align}

From these formulas, the following proposition is evident.

\begin{prop}\label{prop:9.3}
For any knot $K(r)$ in $H(p)$,  
$\mu(r)  \equiv -1$ {\rm (mod  4)}.
\end{prop}

\begin{ex}
(1) Let $p=3$ and $r=[3,4,3,2,3]$.  
Then $\mu(r)= -41$.

(2) Let $p=5$ and $n=2$, then $K(19/85)$ 
belongs to $H(5)$. 

Let $s_0$ be a root of  $1 + 3z + z^2 = 0$.  
Then $b_2  = 2 + s_0$.   
Since $19/85 = [5,2,10]$, it follows
from (9.2)(2) $\mu(r)=8(2+s_0)-1=8s_0+15$
and    
$\lambda_{\rho,K(19/85)}(-1) = (8s_0 + 15)^2$.  

 (3)  Let $p=7$ and $n=3$, 
 then $K(29/217)$ belongs to $H(7)$. 
 Let $s_0$   
be a root of  $1+ 6z +5z^2+ z^3 = 0$.  
Then $b_3 = 3 + 4 s_0 + s_0^2$.   
Since $29/217 = [7,-2,14]$, we have from (9.2)(2)
$\mu(r)=-8(3+4s_0+s_0^2)-1=-25-32s_0-
8s_0^2$ and 
$\lambda_{\rho,K(29/217)}(-1) 
= (25+32s_0+8s_0^2)^2$.
\end{ex}

Two continued fractions $r = [pk_1, 2m_1, pk_2, 2m_2, \cdots, pk_{{\ell}+1}]$ and  $r' = [pk'_1, 2m'_1, pk'_2, 2m'_2, 
\cdots, pk'_{q+1}]$ are said to be (mod 4)-{\it equivalent} 
if $r$ is transformed into $r'$ by a finite sequence of the following 
four operations and their inverses: \\
(1)  replacement of  $k_i$  by $k_i +4d$, $d \in  \ZZ$,\\
(2)	reduction of $[ \cdots, pk_i, 0, pk_{i+1}, \cdots ]$ to $[\cdots, p(k_i + k_{i+1}), \cdots ]$,\\
(3)	reduction of  $[\cdots, 2m_i, 0, 2m_{i+1}, \cdots ]$ to $[\cdots, 2(m_i + m_{i+1}), \cdots ]$,\\
(4) reduction of $[\cdots, pk_r, 2m_r, 0]$ to
$[\cdots, pk_r]$\\
(5) reduction of $[0,2m_1,\dots]$ to $[pk_2,2m_2,\dots]$.

For example, $[p,2,4p,-2,2p]$ is equivalent to $[3p]$.

Computations show that the following conjecture is plausible.

\begin{yosou}
$\mu(r) = -1$ if and only if $r$ is 
(mod 4)-equivalent to either $[p]$ or $[3p]$.
\end{yosou}

\section{Generalization and Silver-Williams Conjecture}

Let $\rho : G(K(r)) \rightarrow SL(2, \ZZ[s_r]) 
\subset SL(2, \CC)$  be a canonical parabolic 
representation of $G(K(r))$ defined in Section 2, 
where $r = \beta / \alpha  <1$. 
The representation polynomial $a(z)$ of $\rho$
has the following properties. (See \cite{R2}.)

\begin{align}
&(1)\  a(z)\ {\rm is\  a\  monic\  integer\ 
 polynomial\  of\  
degree}\  (\alpha - 1)/2.\nonumber\\   
&(2)\	{\rm All\ the\ roots\ of}\ 
a(z) = 0\ {\rm are\ distinct\ and\ simple}.
\end{align}

Let $\widetilde{\Delta}_{\rho,K(r)}(t)$ 
be the twisted Alexander polynomial of $K(r)$ 
associated to  $\rho$. 
Then $\widetilde{\Delta}_{\rho,K(r)}(t)$ is a polynomial 
over $\ZZ[s_r]$.
In order to emphasize this fact, sometimes
we denote it by $\widetilde{\Delta}_{\rho,K(r)}(t|s_r)$.
Let $\theta(z)$ be the minimal polynomial of $s_r$
and deg$\theta(z)=d$.
Let $\gamma_1, \gamma_2, \cdots, 
\gamma_{d}$ be all the roots of 
$\theta(z) = 0$. 
Recently, D.Silver and S.Williams consider the 
integer polynomial $D_{\rho(\theta),K(r)}(t)$ defined as
\begin{equation}
D_{\rho(\theta),K(r)}(t) = \prod_{j=1}^{d} 
\widetilde{\Delta}_{\rho,K(r)}(t|\gamma_j).
\end{equation}

They call it {\it the total $\rho(\theta)$-twisted Alexander 
polynomial} of $K$
and they propose 
the following conjecture.

\begin{yosou} \cite{SW}\label{conj:10.1}
For any 2-bridge knot $K(r)$ and a 
canonical parabolic representation $\rho$, \\
(1) $|D_{\rho(\theta),K(r)}(1)|=2^d$ and\\
(2) $|D_{\rho(\theta),K(r)}(-1)|=2^d N^2$,
where $d={\rm deg}\theta$ and $N$ is a non-zero integer.
\end{yosou}

As they point out, $D_{\rho(\theta),K(r)}(t)$ 
can be evaluated as follows.\\
Let $C$ be the companion matrix of the 
polynomial $\theta(z)$ and 
consider the homomorphism
$\Psi : \ZZ G(K(r)) \rightarrow  
M_{2d, 2d}(\ZZ[t^{\pm 1}])$,
defined by 
$\Psi: x \mapsto 
\mtx{E}{E}{0}{E}$,
$y \mapsto 
\mtx{E}{0}{C}{E}$,
where $E$ is the identity matrix of degree 
$d$.   \\
It is known that 
\begin{equation}
D_{\rho(\theta),K(r)}(t)=
\det [ \widetilde{\Delta}_{\rho,K(r)}(t|C)],
\end{equation}
where $\widetilde{\Delta}_{\rho,K(r)}(t|C)$ 
is a matrix of degree $2d$ obtained from 
$\widetilde{\Delta}_{\rho,K(r)}(t|s_r)$ 
by substituting $C$ for $s_r$.
Computations below show that the conjecture holds for 
$r = 3/5, 3/7$ and $5/9$.  See Example \ref{ex:2.1}.

For $r=3/5$, 
$D_{\rho(\theta),K(r)}(t)=(1-4t+t^2)^2$ 
and hence 
$D_{\rho(\theta),K(r)}(1)=2^2$ and
$D_{\rho(\theta),K(r)}(-1)=2^2 3^2$.
Note that $\theta(z)=a(z)=1-z+z^2$.

For $r=3/7, \widetilde{\Delta}_{\rho,K(r)}(t)
= - (4+s_r^2) + 4t- (4+s_r^2)t^2$, 
and hence, we have 
$D_{\rho(\theta),K(r)}(t)
=\det[\widetilde{\Delta}_{\rho,K(r)}(t|C)]
=25 - 104t + 219t^2 -272t^3 + 
219t^4 - 104t^5 + 25t^6$, 
and  
$D_{\rho(\theta),K(r)}(1)=2^3$
and $D_{\rho(\theta),K(r)}(-1)=2^3 11^2$.
Note $\theta(z)=a(z)=1+2z+z^2+z^3$.

For $r= 5/9, D_{\rho(\theta),K(r)}(t)
= 41 - 376t + 1428t^2 -2984t^3 + 
3798t^4 - 2984t^5 + 1428t^6 - 376t^7 +41t^8$, 
and hence, 
$D_{\rho(\theta),K(r)}(1)=2^4$
and $D_{\rho(\theta),K(r)}(-1)=2^4 29^2$.
Note deg$\theta=4$.

In this section, as a simple application of our main theorem, 
we prove Conjecture \ref{conj:10.1} for a torus knot $K(1/p)$
and a knot $K(r)$ in $H(p)$.

Let $\tau: G(K(1/p))\rightarrow
SL(2,\ZZ[s_0])\subset SL(2,\CC)$ be
the canonical parabolic presentation,
and $a_n(z)$ the representation polynomial of $\tau$.
The properties of $a_n(z)$ are well-studied in 
\cite{R2} and \cite{S}, 
some of which are listed below.

\begin{prop}\label{prop:10.2}
Let $p= 2n+1$.
(1)  $a_n(z) = \prod \chi_s (z)$,
where the product runs over all odd integers 
$s$ dividing $p, 3 \leq s \leq p$ and 
$\chi_s(z)$ is an irreducible, monic integer polynomial. 
The degree of  $\chi_s(z)$ is given by 
$\phi (s)/2$, where $\phi (s)$ 
is Euler function, 
i.e. the number of integers $m, 1 \leq  m \leq s$,
that are relatively prime to $s$. 
In particular, if $p$ is prime, 
then  $\chi_p (z) = a_n(z)$.\\
(2)  $a_n(z) = \sum_{k=0}^n \binom{n+k}{2k}z^k$. 
(3) All the roots of $a_n(z) =0$ are distinct and simple, 
and they are
$- 4 sin^2 \frac{(2k-1)\pi}{2(2k+1)},
1 \leq k \leq n$,
and hence all the roots are real and are in the 
interval $(-4,0)$.
\end{prop}

\begin{ex}
Here are some examples of $\chi_s(z)$.\\
(1)   $\chi_3(z) = a_1(z) = 1+z$\\
(2)   $\chi_5(z) = a_2(z) = 1 +3 z  +z^2$\\
(3)  $\chi_7(z) = a_3(z) = 1 + 6 z  + 5z^2 + z^3$\\
(4)  $\chi_9 (z) = 1 + 9 z  + 6 z^2 + z^3$\\
(5)   $\chi_{15} (z) = 1 + 24 z  + 26 z^2 + 9 z^3  + z^4$\\
(6) $\chi_{21}(z) = 1 + 48 z  + 148 z^2 + 146 z^3 + 
64 z^4 + 13 z^5  +z^6$
\end{ex}

Now let $s_0$ be a zero of $\chi_q(z),q|p,q\ge 3$.
Let $r_1,r_2,\dots,r_d,d={\rm deg}\chi_q(z)=
\phi(q)/2$, be the roots of $\chi_q(z)=0$.
Then, by Proposition \ref{prop:2.1}, 
the total $\tau(\chi_q)$-twisted Alexander polynomial
 $D_{\tau(\chi_q),K(r)}(t)$
 is given by\\
$D_{\tau(\chi_q),K(1/p)}(t)\\
\hspace*{5mm}
=\prod_{j=1}^d[b_1(r_j)+b_2(r_j)t^2+\cdots
+b_n(r_j)t^{2n-2}+b_n(r_j)t^{2n}+\cdots
+b_1(r_j)t^{4n-2}]$,
and hence, by (4.3)(2), we have,\\
$D_{\tau(\chi_q),K(1/p)}(\pm 1)\\
\hspace*{5mm}
=\prod_{j=1}^d [b_1(r_j)+b_2(r_j)+\cdots
+b_n(r_j)+b_n(r_j)+ \cdots
+b_1(r_j)]=\prod_{j=1}^d(-2r_j^{-1})$.

Since deg$\chi_q(z)=d$, we have $r_1 r_2 \cdots r_d=(-1)^d$
and hence\\
 $D_{\tau(\chi_q),K(1/p)}(\pm 1)=2^d$.
This proves Conjecture \ref{conj:10.1} for $K(1/p)$.

Similar arguments work for $K(r)$ in $H(p)$.

Let $\rho=\tau\varphi$ be the canonical parabolic
presentation of $G(K(r))$,\\
$\rho: G(K(r))\rightarrow G(K(1/p))
\rightarrow SL(2,\ZZ[s_0])$.\\
As before, we assume that $s_0$ is a zero of $\chi_q (z),
q|p$, and $r_j,1 \le j\le d$,
are roots of $\chi_q(z)=0$. Then
$\widetilde{\Delta}_{\rho,K(r)}(t|s_0)
=\lambda_{\rho,K(r)}(t|s_0)
\widetilde{\Delta}_{\tau,K(1/p)}(t|s_0)$, and
$D_{\rho(\chi_q),K(r)}(t)
=\{\prod_{j=1}^d \lambda_{\rho,K(r)}(t|r_j)\}
D_{\tau(\chi_q),K(1/p)}(t)$.\\
Now by Theorem A, Propositions \ref{prop:2.1} 
and \ref{prop:4.2}(III)(3), we have\\
$D_{\rho(\chi_q),K(r)}(1)=D_{\tau(\chi_q),K(1/p)}(1)=2^d$.\\
Further, if we write $\lambda_{\rho,K(r)}(-1|r_j)=
\mu_j^2$, then
$D_{\rho(\chi_q),K(r)}(-1)=(\mu_1\mu_2\cdots\mu_d)^2 2^d$.
This proves Conjecture \ref{conj:10.1} for the
total $\rho(\chi_q)$-twisted Alexander polynomial of $K(r)$ in $H(p)$.

\begin{prop}\label{newProp:10.4}
For a knot $K(1/p)$, the total $\rho(\chi_q)$-twisted Alexander polynomial
$D_{\rho(\chi_q),K(r)}(t)$ can be determined 
by the following three formulas.
%
Let $p=2n+1$.\\
(1) If $q$ is a divisor of $p$, say $p=vq, v\ge 3$,
then \\
$D_{\tau(\chi_q),K(1/p)}(t)=
(1-t^{2q}+t^{4q}-\cdots+t^{2(v-1)q})^{d_q}
D_{\tau(\chi_q),K(1/p)}(t)$, where $d_q=$deg$\chi_q(t)$.\\
(2) $\prod D_{\tau(\chi_u),K(1/p)}(t)=(1+t^2)(1+t^{4n+2})^{n-1}$,
\hfill (10.4)\\
where the product runs over all divisors $u(\neq1)$ of $p$.\\
(3) If $p$ is a prime, then 
$D_{\tau(\chi_p),K(1/p)}(t)=(1+t^2)(1+t^{4n+2})^{n-1}$.
\end{prop}

Since 
Proposition \ref{newProp:10.4}(1)
is an easy consequence of Proposition
\ref{prop:3.4} 
and Proposition \ref{newProp:10.4}(3) follows from Proposition \ref{newProp:10.4}(2), only a proof of
Proposition \ref{newProp:10.4}(2) will be given in Appendix (III).

\begin{ex}\label{ex:10.4}
(1) Let $p=9$ and $n=4$.
Then $a_4 (z)= \chi_3 (z) \chi_9(z)$. 
First, by Proposition \ref{newProp:10.4}(3)
$D_{\tau(\chi_3), K(1/3)}(t)=1+t^2$, and by 
Proposition \ref{newProp:10.4}(1),
we see\\
$D_{\tau(\chi_3), K(1/9)}(t)=(1+t^2)(1-t^6 + t^{12})$. 
Further, by 
Proposition \ref{newProp:10.4}(2),\\
$D_{\tau(\chi_3), K(1/9)}(t) D_{\tau(\chi_9), K(1/9)}(t)=(1+t^2)(1+t^{18})^3$, and
hence,\\
$D_{\tau(\chi_9), K(1/9)}(t)=(1+t^2)(1+t^{18})^3 /(1+t^2)(1-t^6 +t^{12})
 = (1+t^{18})^2 (1+t^6)$.\\
%
(2) Let $p=15$ and $n=7$.
Then $a_7(z)=\chi_3 (z) \chi_5(z) \chi_{15}(z)$
and
\begin{align*}
D_{\tau(\chi_3), K(1/15)}(t)
&=\lambda_{\tau(\chi_3), K(1/15)}(t) D_{\tau(\chi_3), K(1/3)}(t)\\
&= (1+t^2)(1-t^6+t^{12}-t^{18}+t^{24}),
\end{align*}
\begin{align*}
D_{\tau(\chi_5), K(1/15)}(t)
&=\{\lambda_{\tau(\chi_5), K(1/15)}(t)\}^2 D_{\tau(\chi_5),K(1/5)}(t)\\
&=(1-t^{10} +t^{20})^2 (1+t^2)(1+t^{10}).
\end{align*}
Since $\prod_{j=3,5,15}D_{\tau(\chi_j), K(1/15)}(t)=(1+t^2)(1+t^{30})^6$, we have
\begin{align*}
D_{\tau(\chi_{15}), K(1/15)}(t)
&=(1+t^2) (1+t^{30})^6 / D_{\tau(\chi_3), K(1/15)}(t)
D_{\tau(\chi_5), K(1/15)}(t)\\
&=(1-t^2+t^4)(1+t^{10})(1+t^{30})^3.
\end{align*}
\end{ex}

Next we discuss a generalization of our main theorem.\\
Suppose there is an epimorphism 
$\varphi$ from $G(K(r))$ to $G(K(r_0))$. 
Using a canonical parabolic representation $\rho$
of $G(K(r_0))$ into $SL(2,\CC)$, 
we can define the
the twisted Alexander polynomials  
$\widetilde \Delta_{\rho \varphi, K(r)}(t)$ 
and   
$\widetilde \Delta_{\rho , K(r_0)}(t)$               
associated to $\rho \varphi$ and $\rho$, 
respectively.  
Since  $\widetilde \Delta_{\rho , K(r_0)}(t)$              
divides  
$\widetilde \Delta_{\rho \varphi, K(r)}(t)$, 
the quotient 
$\lambda_{\rho , K(r)}(t)$ is well-defined.
The following  conjecture is a generalization 
of our main theorem.   

\begin{yosou}  
(1)  $\lambda_{\rho , K(r)}(1) = 1$, 
and\\
(2)  $\lambda_{\rho , K(r)}(-1) = \mu^2$  
for some  $\mu \in \ZZ[s_{r_0}]$.
\end{yosou}

In fact, there is an epimorphism 
$\varphi : G(K(63/115))  \rightarrow G(K(3/5))$ 
and we have the twisted Alexander polynomial of 
$K(63/115)$\\
$\widetilde \Delta_{\rho \varphi, K(63/115)}(t) = 
\lambda_{\rho, K(63/115)}(t) 
\widetilde \Delta_{\rho, K(3/5)}(t)$,
where   $\lambda_{\rho, K(63/115)}(t) 
=  (3 - w) - (16 - 8w) t   + (33 - 34w) t^2    
- (40 - 76w) t^3  + (41 - 98w) t^4   
- (40 - 76w) t^5   + (33 - 34w) t^6    
 - (16 - 8w) t^7   + (3 - w) t^8$, 
and $w$ is a primitive third root of $1$.  
We see then 
$\lambda_{\rho, K(63/115)}(1)=1$ and
$\lambda_{\rho, K(63/115)} (-1)= 225-336 w 
= (17-8w)^2$.

Finally, we give a few remarks on the representation 
polynomials.
Let $f(z)$ and $g(z)$, respectively, be the representation polynomials of $\rho\varphi$
and $\rho$.
Then $g(z)$ divides $f(z)$.
However, the converse seems quite likely to hold,
and therefore, we propose the following conjecture.

\begin{yosou}\label{conj:10.6}
Let $f_1(z)$ and $f_2(z)$, respectively, be the representation polynomials of the canonical parabolic representations
$\rho_1${\rm :} $G(K(r_1))\rightarrow SL(2,\CC)$
and $\rho_2${\rm :} $G(K(r_2))\rightarrow SL(2,\CC)$.
If $f_2(z)$ divides $f_1(z)$, then there exists an epimorphism from
$G(K(r_1))$ to $G(K(r_2))$.
\end{yosou}

It is proven \cite{Kn},\cite{ASWY} that Conjecture
\ref{conj:10.6} holds if $r_2=1/p$
or equivalently, if $\chi_p(z)$ divides $f_1(z)$,
then there exists an epimorphism
from $G(K(r_1))$ to $G(K(1/p))$.

\begin{rem}
Very recently we learned \cite{SW2}
that D. Silver and S. Williams
proved Conjecture \ref{conj:10.1} (1) for any $2$-bridge knot 
$K(r)$.
\end{rem}

\section{Appendix}

In Appendix, we discuss four topics.

(I) Outline of the proof of Proposition \ref{prop:2.1}.\\ 
Consider a Wirtinger presentation of $G(K(1/p))$ given by (2.2):

$G(K(1/p)) = \langle x, y | R_0 = W x W^{-1}y^{-1}\rangle$, 
where $p=2n+1$ and $W = (xy)^n$.  
Then       
$\frac{\partial R_0}{\partial x} 
= (1-y) Q_{n-1} + (xy)^n$ and hence,
\begin{align*}
D&= 
(\frac{\partial R_0}{\partial x})^{\Phi}\\
&=\mtx{1-t}{0}{-s_0 t}{1-t}
\mtx
{\sum_{k=0}^{n-1}a_k t^{2k}  } 
{\sum_{k=0}^{n-1} b_k t^{2k}}
{\sum_{k=0}^{n-1} c_k t^{2k}}
{\sum_{k=0}^{n-1}d_k t^{2k}}
+\mtx{0}{b_n t^{2n}}{c_n t^{2n} }{d_n t^{2n}}\\
&=\mtx{h_{11}}{h_{12}}{h_{21}}{h_{22}},\ {\rm where}
\end{align*}
\begin{align*}
&h_{11}=(1-t) \sum_{k=0}^{n-1} a_k t^{2k},\\
&h_{12}=(1-t) \sum_{k=0}^{n-1} b_k t^{2k} + b_n t^{2n},\\
&h_{21}=-s_0 t \sum_{k=0}^{n-1} a_k t^{2k} 
+ (1-t) \sum_{k=0}^{n-1}c_k t^{2k} +c_n t^{2n}, \ {\rm and}\\
&h_{22}
=
-s_0 t \sum_{k=0}^{n-1}b_k t^{2k} 
+(1-t) \sum_{k=0}^{n-1}d_k t^{2k} + d_n t^{2n}.
\end{align*}

Since for 
$j \geq 1, a_0 + a_1 + \cdots + a_{j- 1} = b_j$ and 
$s_0 b_j = c_j$, we see that 
\begin{equation*}
 \sum_{k=0}^{n - 1} a_k t^{2k} = \sum_{k=1}^n b_k t^{2k - 2}(1 - t^2) + b_n  t^{2n},
 \ {\rm and\ hence}
 \end{equation*} 
\begin{equation*}
h_{2,1} = (1-t) \{ -  \sum_{k=1}^n c_k t^{2k - 1} (1+t) 
+ \sum_{k=0}^{n - 1}c_k t^{2k} + c_n t^{2n}\}
= - (1-t) \sum_{k=1}^n c_k t^{2k - 1}.
\end{equation*}      

Thus the first column is divisible by $1-t$ and hence,
\begin{align*}
&\det D/(1-t)
=\det\mtx{h'_{1,1}}{h_{1.2}}{h'_{2,1}}{h_{2,2}}, \ {\rm where}\\  
&h'_{1,1}
=\sum_{k=0}^{n - 1} a_k t^{2k}\ {\rm and}\  
h'_{2,1}
=- \sum_{k=1}^n c_k t^{2k - 1}.   
\end{align*}

Now subtract the first column multiplied through 
$t$ from the second column so that we have
\begin{align*}
h_{1,2} - t h'_{1,1}
&=b_n t^{2n} + (1-t) \sum_{k=1}^{n - 1}b_k t^{2k} 
- \sum_{k=1}^n b_k t^{2k - 1} (1- t^2) 
- b_n t^{2n+1}\\
&=- (1-t) \sum_{k=1}^n b_k t^{2k - 1}.
\end{align*}    

Similarly, noting 
$c_k + d_k = a_k$ for  $k \geq 0$, we obtain
\begin{align*}
h_{2,2} -  t h'_{2,1}
&=- t \sum_{k=0}^{n - 1} c_k t^{2k}    
+ (1-t) \sum_{k=0}^{n - 1}d_k t^{2k} + d_n t^{2k} 
+ t \sum_{k=1}^n c_k t^{2k - 1}\\
&=(1-t) \sum_{k=0}^{n - 1}a_k t^{2k}.
\end{align*} 

Therefore
\begin{align*}
\det D/(1-t)^2 
=
\det 
\mtx{\sum_{k=0}^{n-1}a_k t^{2k}}
{-\sum_{k=1}^n b_k t^{2k-1}  }
{ - \sum_{k=1}^n c_k t^{2k-1}}
{ \sum_{k=0}^{n-1}a_k t^{2k} },\ {\rm and\ hence}\\
\widetilde \Delta_{\rho, K(1/p)}(t)
=
( \sum_{k=0}^{n-1}a_k t^{2k})^2 -  
( \sum_{k=1}^n b_k t^{2k-1}) 
( \sum_{k=1}^n c_k t^{2k-1}).
\end{align*}     
Proposition \ref{prop:2.1}, 
then, follows from (A.1) below:

\hfill
For $k=1, 2, \dots, n$,
$b_k = \sum_{i+j=k-1} a_i a_j - \sum_{i+j=k}b_i c_j$.\hfill(A.1)  \\   

(II) 
Proof of Proposition \ref{prop:3.4}.\\ 
We prove

\hfill
$\lambda_{\rho, K(1/pq)}(t)
=\Delta_{K(1/q)} (t^{2p})$.\hfill (A.2) 

Let $p=2n+1$ and $q=2m+1$. 
Let
$G(K(1/pq))= \langle x,y| R_{pq}\rangle$ 
and  $G(K(1/p)) = \langle x,y| R_0\rangle$
be Wirtinger presentations,
where  $R_{pq} =
(xy)^{2mn+m+n} x (xy)^{-(2mn+m+n)}y^{-1}$  and
$R_0=(xy)^n x (xy)^{-n}y^{-1}$.
We must express $R_{pq}$ as a product of conjugates of 
$R_0$.  In fact,
we prove:

\medskip

\noindent
{\bf Lemma A.1.}  
{\it
$R_{pq}=R_0^{\tau_m}$, where 

\hfill
  $\tau_m=
\sum_{k=0}^m (xy)^{kp} 
-  \sum_{k=0}^{m-1} (xy)^{kp}(xy)^n x$.
\hfill(A.3)  
}

\medskip
\noindent
{\it Proof.}  
We prove (A.3) by induction on $m$.
If $m = 1$, then 
$R_{pq} =(xy)^{3n+1} x (xy)^{-(3n+1)}y^{-1}
=
R_0^{\tau_1}$, where
$\tau_1=(xy)^{2n+1} - (xy)^n x  + 1$.
Thus (A.3) holds.
Now inductively, consider  $\tau_{m+1}$.  
Applying the previous 
argument repeatedly, we obtain 
\begin{align*}
R_{p(2m+3)}
&=(xy)^{mp+3n+1}x (xy)^{-(mp+3n+1)}y^{-1}\\
&=(xy)^{mp+p} \{(xy)^n x (y^{-1}x^{-1})^n y^{-1}\}
x^{-1}(y^{-1}x^{-1})^{mp+2n}y^{-1}\\
&=R_0^{(xy)^{(m+1)p}}(xy)^{(m+1)p} 
x^{-1}(y^{-1}x^{-1})^{mp+2n}y^{-1}\\      
&=R_0^{(xy)^{(m+1)p}}(xy)^{mp+n} x 
\{y (xy)^n x^{-1}(y^{-1}x^{-1})^n \}
(y^{-1}x^{-1})^{mp+n} y^{-1}\\
&=R_0^{(xy)^{(m+1)p}} R_0^ {-(xy)^{mp+n} x}
(xy)^{mp+n} x (y^{-1} x^{-1})^{mp+n} y^{-1}\\
&=R_0^{(xy)^{(m+1)p} -(xy)^{mp+n} x} R_0^{\tau_m},\        
{\rm and\  hence}
\end{align*}
$\tau_{m+1}=(xy)^{(m+1)p} - (xy)^{mp}(x y)^n x + \tau_m$.
This proves (A.3).
\qed

Now to evaluate  $\lambda_{\rho,K(1/pq)}(t)$, 
we compute $ \Phi (\tau_m)$ that is given 
as follows.   
Since  
$\Phi ((xy)^{kp})=\mtx{(-1)^k}{0}{0}{(-1)^k} t^{2kp}$
and   
$\Phi [(xy)^n x] = 
\mtx{0}{b_n}{c_n}{0} t^p$, we have
\begin{align*}
\Phi (\tau_m)  
&= \sum_{k=0}^m 
\mtx{(-1)^k}{0}{0}{(-1)^k} t^{2kp}-
 \sum_{k=0}^{m-1}
 \mtx{0}{(-1)^k b_n}{(-1)^k c_n}{0}t^{(2k+1)p} \\
&  = 
\left[
\begin{array}{ll}
\ \ \ {\displaystyle \sum_{k=0}^m} (-1)^k t^{2kp}
&\ \ 
    -  {\displaystyle \sum_{k=0}^{m-1}} (-1)^k b_n t^{(2k+1)p}\\
     - {\displaystyle \sum_{k=0}^{m-1}} (-1)^k c_n t^{(2k+1)p}
&\ \ \ \  {\displaystyle \sum_{k=0}^m} (-1)^k t^{2kp}
\end{array}
\right].
    \end{align*}
Since $b_n c_n = -1$, we see

\begin{align*}
\det
[\Phi (\tau_m)]
&=
\left\{ \sum_{k=0}^m (-1)^k t^{2kp}]^2  + \sum_{k=0}^{m-1} 
(-1)^k t^{(2k+1)p}\right\}^2\\
&=\sum_{k=0}^{2m}(-1)^k t^{2kp}\\
&=\Delta_{K(1/q)} (t^{2p}).
\end{align*}

This proves (A.2).\\

(III) Sketch of the proof of Proposition \ref{newProp:10.4}(2).

Denote $x_j = t^{2j} + t^{4n - 2j - 2}$, 
$0 \leq j \leq n-1$.   Then 
$\widetilde{\Delta}_{\rho,K(1/p)}(t)$ can be written as 
$\widetilde{\Delta}_{\rho,K(1/p)}(x_0, \cdots, x_{n - 1})
=b_1 x_0  +  b_2 x_1 + b_3 x_2 + \cdots + b_n x_{n - 1}$.\\
We use the following easy formula proved in \cite{S}.

\hfill
 For $k \geq 1,
b_k=\sum_{j=0}^{k - 1} \binom{k+j}{2j+1}s_0^j$,   
where $s_0$ is a root of $a_n (z)$.\hfill (A.3)  \\
Let $C  = \left[ c(i,j)\right]_{1 \leq i,j \leq n}$ 
be the companion matrix of  $a_n (z)$.          \\
Only non-zero entries of $C$ are: \\
\hspace*{2cm}(1)	For $i = 1,2, \cdots, n,
c(i, n) = - \binom{n+i - 1}{2(i - 1)}$\\     
\hspace*{2cm}(2) For $1 \leq i \leq  n-1,c(i+1,i) = 1$. 
\hfill \mbox{(A.4)}\\
%
Let $C^k  = \left[ c_k (i,j)\right]_{1 \leq i,j \leq n}$.  
Then a straightforward calculation verifies 
the following lemma.

\medskip
\noindent
{\bf Lemma A.2.  } 
{\it Let $1 \leq k \leq n-1$.\\
(1)   For  $1 \leq i \leq n-k ,  c_k (k+i,i) = 1$.\\
(2)  For $1 \leq i \leq n,  1 \leq k \leq n-1,  
1\leq j \leq  k-1,
c_k (i, n-k+j) = 
c_{k - 1} (i, n- k+j+1) = 
c_{k - 2}(i,n-k+j+2) =  \cdots 
=  c_{j+1}(i,n-1) = c_j (i,n)$.\\
(3) For $k \geq 2, c_k (1,n) = 
c_1(1,n) c_{k - 1}(n,n)$, and 
for $i \geq 2$ and  $k \geq 2$,
$c_k (i,n) = c_{k - 1}(i-1,n) 
+ c_1(i,n) c_{k - 1}(n,n)$.\\
(4) Other values of  $c_k (i,j)$ are $0$.
}

\medskip
Let ${\displaystyle B_k = \sum_{j=0}^{k - 1} 
\binom{k+j}{2j+1}C^j}$ and
${\displaystyle D = \sum_{j=1}^n B_j x_{j - 1} = 
\left[ d(i,j)\right]_{1 \leq i,j \leq n}
}$.
Since $a_n(z)$ is separable, det$D$ is the LHS of 
(10.4).
We determine $d(i,j)$.
Since the following three lemmas are easily proven, 
we omit the details.

\medskip
\noindent
{\bf
Lemma A.3.}{\it
(1)   For  $1 \leq i < j \leq n$,
\begin{equation*}
d(i,j) =  \sum_{m=1}^{j - 1}
\sum_{k=0}^{j - 2 - m+1}
\binom{2n - 2j + 2m+1+k}{2n - 2j+2m+1} 
c_m (i,n) x_{n - j+m+k},
\end{equation*}
(2)   For $1 \leq i \leq n$,
\begin{equation*}
d(i,i) = \sum_{k=1}^n k x_{k - 1} 
+ \sum_{m=1}^{i - 1} \sum_{k=0}^{i - m - 1} 
\binom{2n - 2i+2m+1+k}{2n - 2i+2m+1} 
c_m (i,n) x_{n - i+m+k},
\end{equation*}
(3)   For $1\leq j < i \leq n$,
\begin{align*}
d(i,j) &= \sum_{m=1}^{j - 1} 
\sum_{k=0}^{i - m - 1} 
\binom{2n - 2j+2m+1+k}{2n - 2j+2m+1} 
c_m (i,n) x_{n - j+m+k}\\
 & \ \ + 
\sum_{k=0}^{n - i+j-1}
\binom{2i - 2j+1+k}{2i - 2j+1}x_{i - j+k}.
\end{align*}
(4)  For $1\leq j \leq n-1$,
\begin{align*}
d(n,j)&=\sum_{m=1}^{j - 1} 
\sum_{k=0}^{j - m - 1} 
\binom{2n - 2j+2m+1+k}{2n - 2j+2m+1} 
c_m (i,n) x_{n - j+m+k} \\
& \ \ + \sum_{k=0}^{j-1}
\binom{2n - 2j+1+k}{2n - 2j+1}x_{n - j+k}.
\end{align*}
In particular, $d(n,1) = x _{n - 1}$.    
}
\medskip

There are some relations among entries of $D$.

\noindent{\bf Lemma A.4.} {\it
For  $2 \leq i,j \leq n$\\
(1) $d(1,j) = c_1(1,n) d(n,j-1)$,\\
(2) $d(i,j) = d(i-1,j-1) + c_1(i,n) d(n,j-1)$.                                          
}\\

Note that for 
$k \ne n-1, x_k t^{-2}(1-t^2)^2 
= x_{k - 1} - 2 x _k + x_{k+1}$ 
and $x_{n - 1}t^{-2}(1-t^2)^2 
= x_{n - 2} - x_{n - 1}$.    

\medskip
Using this lemma, we can prove:

\noindent
{\bf
Lemma A.5.}
{\it
$d(i,i) t^{-2}(1-t^2)^2
=d(i, i+1) + t^{-2}(1+t^{4n+2})$, and 
if $j \ne  i$, 
then $d(i,j) t^{-2}(1-t^2)^2 = d(i,j+1)$. 
}

\medskip
Now consider 
$D= \left[ d(i,j) \right]_{1 \leq i, j \leq n}$.
First subtract the $(n-1)^{\rm st}$ column multiplied
through $t^{-2}(1-t^2)^2$ 
from the $n^{\rm th}$ column. 
Then by Lemma A.4, all the entries of the 
resulting $n^{\rm th}$ column are $0$ 
except the $(n-1, n)$ entry that is 
$-  t^{-2}(1+t^{4n+2})$.  
Successive applications of the same operation 
applied on the $(j-1)^{\rm st}$ column and the 
$j^{\rm th}$ column transform 
$D$ into a new matrix 
$\widehat{D} 
= \left[ \widehat{d}(i,j) \right]_{1 \leq i, j \leq n}$, 
where the off diagonal entries 
$\widehat{d}(i,i+1)$ are 
$t^{-2}(1+t^{4n+2})$ and 
$\widehat{d}(i,n) = d(i,n), 1 \leq i  \leq n$ 
and all the rest is $0$.
Thus\begin{align*}
\det D
&=\det \widehat{D}\\
&=\{- t^{-2}(1+t^{4n+2})\}^{n - 1} 
(-1)^{n - 1}x_{n - 1}\\
&=t^{-(2n - 2)}(1+t^{4n+2})^{n - 1}t^{2n-2}(1+t^2)\\
&=(1+t^2) (1+t^{4n+2})^{n - 1}.
\end{align*}

%

\medskip
(IV) Alternative characterization of
$r$ for $K(r)$ in $H(p)$.

\noindent{\bf Definition A.6.}
Let $\alpha$
and $\beta$ be co-prime odd integers
with $0<|\beta|<\alpha$, and $p$ an odd integer.\\
(I)
We say that $r=\beta/\alpha$ is 
$p$-{\it expandable} if $r$ has a continued
fraction expansion of the form:\\
\centerline{
$r=[p k_1, 2 m_1,
p k_2, 2m_2, \ldots]$, where $k_i, m_i
\in \ZZ\setminus\{0\}$.}\\
(II)
We know that $r$ has a unique continued 
fraction\\
\centerline{
$r=[2a_1, 2a_2, \ldots,
2a_\ell, c]$, where $c\neq\pm1$.}
Then we inductively define
$r$ to be  
$p$-{\it admissible} by the following:

(a) $[c]$ is $p$-admissible 
if and only if $c\equiv p \mod 2p$,

(b) $[2a_1, c]$ is never $p$-admissible, and

(c) Let $r=[2a_1, 2a_2, x,\ldots]$, where 
$x, \ldots$ denotes
$2a_3, \ldots$ or $c$.\\
\hspace*{0.9cm}
Then $r$ is $p$-admissible if and only if
one of the following is satisfied:\\
\hspace*{1.5cm}(i) $2a_1\equiv 0 \mod 2p$ and
$[x,\ldots]$ is
$p$-admissible.\\
\hspace*{1.5cm}(ii) $2a_1\equiv p+1 \mod 2p, 2a_2=2$,
and $[x-(p+1), \ldots]$
is $p$-admissible.\\
\hspace*{1.5cm}(iii)  $2a_1\equiv p-1 \mod 2p, 2a_2=-2$,
and $[x-(p-1),  \ldots]$
is $p$-admissible.

\medskip

\noindent{\bf Example A.7.}
Let $r=12225937/33493827$. \\
Then $r$ is both $3$-expandable and 
$3$-admissible, since
\begin{align*}
r&=
[3,4,6,-4,9,6,18,-2,-3,4,6]\\
&
=[2,-2,-2,-2,6,2,2,2,10,6,18,-2,-4,-2,-2,-2,5].
\end{align*}

\medskip

\noindent{\bf Remark A.8.}
(1) Let $p$ be an odd integer.
Then both of the denominator
and numerator of $r=
[pk_1, 2m_1, pk_2, 2m_2, \ldots]$
are odd if and only if 
(i) the length of expansion is odd
and (ii) total of $k_i$'s is odd.
(2) If both of the denominator
and numerator of $r$ is odd, then
the reduction in Definition A.6.  
(c) preserves that property.

\medskip

\noindent{\bf Lemma A.9.}
For continued fractions, we have the following
equalities:\\
(1) $[\dots, a, 2, b, \dots]=[\dots, a-1, -2, b-1,\dots]$\\
(2) $[\dots, a, \underbrace{2, 2, \ldots, 2}_{k},b,\ldots]=[\ldots, a-1, -(k+1), b-1,\ldots]
$

\medskip

\noindent{\bf Theorem A.10.}
{\it
Let $\alpha$
and $\beta$ be co-prime odd integers
with $0<|\beta|<\alpha$, and $p$ an odd integer.
Then $r=\beta/\alpha$ is $p$-admissible
if and only if $r$ is $p$-expandable.
}

\medskip

\noindent
{\it Proof.}
(Proof of \lq$\Rightarrow$\rq)\ 
Suppose $r=[2a_1, 2a_2, \ldots, 2a_\ell, c]$ is
$p$-admissible. 
We prove that $r$ is $p$-expandable
by induction on the length of expansion.
First, if $[c]$ is $p$-admissible,
then $c=2pn +p$ for some $n\in \ZZ$
and hence $r$ is $p$-expandable.
Next, if $r=[2a_1, c]$, $r$ is not $p$-admissible
and there is nothing to prove.
Let $r=[2a_1, 2a_2, x,\ldots]$, where $x$
denotes $2a_3$ or $c$.

Case 1, $2a_1=2pn$ for some $n$:
Here, $[x, \ldots]$
is $p$-admissible, and by 
induction hypothesis,
$[x, \ldots]$ is $p$-expandable. So,
$r=[2pn,2a_2, x, \ldots]$ is also $p$-expandable.

Case 2, $(2a_1,2a_2)=
(2pn+(p+1), 2)$ for some $n$: 
Here, $[x-(p+1), \ldots]$ is $p$-admissible,
and hence by induction hypothesis,
$[x-(p+1), \ldots]$ is $p$-expandable.
Then $[p(2n+1), -2, x-(p+1)+p, \ldots]$ is
also $p$-expandable.
Since 
\begin{align*}
&\ \ \ \ [p(2n+1), -2, x-(p+1)+p,\ldots]\\
&=[p(2n+1)+1, 2, x-(p+1)+p+1,\ldots]\\
&=r,
\end{align*}
we see that $r$ is $p$-expandable.

Case 3, $(2a_1,2a_2)=
(2pn+(p-1), -2)$: 
This case is similar to Case~2.


(Proof of \lq$\Leftarrow$\rq)\ 
Suppose that the length of expansion is
$1$, i.e., $r=[p k_1]$. 
Since $\alpha$ and $\beta$ are odd, 
both $p$ and $k_1$ are odd.
Therefore, writing $k_1=2q+1$,
we see that $p k_1=p(2q+1)\equiv p \mod 2p$,
and hence $r$ is $p$-admissible.
The length of expansion is never equal to
$2$, since if so,
 $r=\frac{1}{pk-1/2m}=\frac{2m}{2pkm-1}$
and hence $\beta$ would be even.\\
Let $r=[pk_1, 2m_1, pk_2,\ldots]$.

Case 1, $k_1$ is even:
Here, we can write $pk_1=2pq$,
and hence
it suffices to show that
$[pk_2, \ldots]$ is 
$p$-admissible, which is true since,
by Remark A.8,
we can use the induction hypothesis.

Case 2.1, $k_1$ is odd and $m_1>0$:
Write $k_1=2q+1$, then we have
\begin{align*}
r
&=[2qp+p, 2m_1, pk_2, \ldots]\\
&=[2qp+p-1, 
\underbrace{-2, \ldots, -2}_{2m_1 -1},\
pk_2-1, \ldots]\\
&=[\{2qp+p-1, -2\}, \underbrace{
\{-2,-2\}, \ldots, \{-2,-2\}}_{m_1-1},\ pk_2-1,\ldots]. 
\end{align*}
(Braces are inserted just for
the sake of pairing.)
Then we further see that
$r$ is $p$-admissible if and only if
so is\\
\centerline{
$[\{-2-(p-1), -2\},
\underbrace{\{-2,-2\},\ldots, \{-2,-2\}}_{
m_1-2}, pk_2-1,\ldots]$.}\\
Since $-2-(p-1) \equiv p-1 \mod 2p$,
we see that $r$ is
$p$-admissible if and only if
so is
$[\{-2-(p-1), -2\},
\underbrace{\{-2,-2\},\ldots, \{-2,-2\}}_{
m_1-3}, pk_2-1,\ldots]$.
Repeatedly
$r$ is
$p$-admissible if and only if
so is
$[-(p-1) +pk_2-1, \ldots]$, which
is $p$-expandable
since $-(p-1)+pk_2-1=p(k_2+1)$.
Now, $k_1+k_2$ and $k_2+1$ have the same
parity and hence, by Remark A.8,
we can use the
induction hypothesis to see that $r$
is $p$-admissible.

Case 2.2,  $k_1$ is odd and $m_1<0$:
This case is similar to Case 2.1.

This completes the proof of Theorem A.10.
\qed

{\bf Acknowledgements. } 
First we would like to express our deep
appreciation to Daniel Silver and Susan Williams 
who give us many invaluable comments 
on our present work. 
Also, we thank Hiroshi Goda who informed us 
Kitayama's work \cite{Ktym} on a refinement
of the invariance of the twisted Alexander 
polynomials of knots.  
Further,
we thank 
Makoto Sakuma for giving us helpful
information regarding this work, and Alexander Stoimenow
who gave us the table of polynomials 
$\lambda_{\rho, K(r)}(t)$
for many $2$-bridge knots in $H(3)$.

The first author is
partially supported by MEXT, Grant-in-Aid for
Young Scientists (B) 18740035,
and the second author is
partially supported by NSERC Grant~A~4034
%


%
%
%
%

\end{document}